\documentclass[reqno]{amsart}
\usepackage[scale=0.75, centering, headheight=14pt]{geometry}
\usepackage[latin1]{inputenc}
\usepackage[T1]{fontenc}
\usepackage{lmodern}
\usepackage[english]{babel}
\usepackage{amsmath}

\usepackage[
    type={CC},
    modifier={by-nc-nd},
    version={4.0},
]{doclicense}
               
\usepackage{amsmath,amssymb,amsfonts,amsthm}
\usepackage{mathtools,accents}
\usepackage{mathrsfs}
\usepackage{xfrac}
\usepackage{array} 
\usepackage{aliascnt}
\usepackage{booktabs} 
\usepackage{array} 

\usepackage{verbatim} 
\usepackage{subfig} 

\usepackage{mathrsfs, dsfont}
\usepackage{amssymb}
\usepackage{amsthm}
\usepackage{amsmath,amsfonts,amssymb,esint}
\usepackage{graphics,color}
\usepackage{enumerate}
\usepackage{mathtools,centernot}
\usepackage{aligned-overset}

\usepackage{microtype}
\usepackage{paralist} 
\usepackage{cases}
\usepackage[initials, alphabetic]{amsrefs}
\allowdisplaybreaks

\usepackage{braket}
\usepackage{bm}

\usepackage[citecolor=blue,colorlinks]{hyperref}
\addto\extrasenglish{}

\usepackage{enumerate}
\usepackage{xcolor}
\usepackage{aliascnt}

\makeatletter
\def\newaliasedtheorem#1[#2]#3{
  \newaliascnt{#1@alt}{#2}
  \newtheorem{#1}[#1@alt]{#3}
  \expandafter\newcommand\csname #1@altname\endcsname{#3}
}
\makeatother

\usepackage{indentfirst}

\usepackage{esint}
\usepackage{mathrsfs}
\usepackage{stmaryrd}

\DeclareUnicodeCharacter{00A0}{ }
\DeclareUnicodeCharacter{00A0}{~}

\makeatletter
\newsavebox{\measure@tikzpicture}

\newcommand{\setword}[2]{%
  \phantomsection
  #1\def\@currentlabel{\unexpanded{#1}}\label{#2}%
}

\renewcommand\labelenumi{(\arabic{enumi})}
\renewcommand\theenumi\labelenumi

\newtheorem{theorem}{\bf Theorem}[section]
\newtheorem{remark}[theorem]{\bf{Remark}}
\newtheorem{definition}[theorem]{\bf Definition}
\newtheorem{lemma}[theorem]{\bf Lemma}
\newtheorem{proposition}[theorem]{\bf Proposition}

\newtheorem*{theorem*}{Theorem}

\newtheorem{maintheorem}{Theorem}

\newcommand{\eps}{\varepsilon}
\newcommand{\R}{\mathbb R}

\newcommand{\rank}[1]{\operatorname{rank}(#1)}

\newcommand{\dist}[2]{\operatorname{dist} ( #1, #2 )}

\newcommand{\diam}[1]{\operatorname{diam} \left(#1 \right)}
\newcommand{\rotational}[1]{\operatorname{curl} \left(#1 \right)}
\newcommand{\divergence}{\operatorname{div}}

\newcommand{\sym}{\operatorname{sym}}
\DeclareMathOperator{\loc}{loc}

\allowdisplaybreaks
\numberwithin{equation}{section}

\title{Wild solutions to scalar Euler-Lagrange equations}

\author[C. J. P. Johansson]{Carl Johan Peter Johansson}

\address{Carl Johan Peter Johansson  
\hfill\break  EPFL SB, Station 8, CH-1015 Lausanne, Switzerland}
\email{carl.johansson@epfl.ch}                                      

\subjclass[2020]{Primary 35D30, 35J60, 35A02; Secondary 35A09} 
\thanks{Author's Accepted Manuscript of an article first published in Transactions of the American Mathematical Society in Volume 377, Number 7, published by American Mathematical Society \copyright2024 American Mathematical Society}

\begin{document}
\doclicenseThis
\maketitle

\begin{abstract}
We study very weak solutions to scalar Euler-Lagrange equations associated with quadratic convex functionals. We investigate whether $W^{1,1}$ solutions are necessarily $W^{1,2}_{\loc}$, which would make the theories by De Giorgi-Nash and Schauder applicable.
We answer this question positively for a suitable class of functionals. This is an extension of Weyl's classical lemma for the Laplace equation to a wider class of equations under stronger regularity assumptions. Conversely, using convex integration, we show that outside this class of functionals, there exist $W^{1,1}$ solutions of locally infinite energy to scalar Euler-Lagrange equations. In addition, we prove an intermediate result which permits the regularity of a $W^{1,1}$ solution to be improved to $W^{1,2}_{\loc}$ under suitable assumptions on the functional and solution.
\end{abstract}

\section{Introduction}\label{sec:Intro}
Let $f \in C^2(\R^n)$ be uniformly convex of quadratic growth and $\Omega \subset \R^n$ any bounded Lipschitz domain. The equation
\begin{equation}\label{eq:Intro:EL}
\tag{EL}
 \divergence (Df(Du)) = 0 \text{ in } \Omega
\end{equation}
is the Euler-Lagrange equation of the energy functional $\mathcal{E}_f \colon W^{1,2}(\Omega) \to \R$ defined by
\[
 \mathcal{E}_f(u) \coloneqq \int_{\Omega} f(Du) \, dx.
\]
This is a standard model in the Calculus of Variations. The classical regularity theory was obtained in the fundamental works by E. De Giorgi in \cite{EDG} and J. Nash in \cite{JN58}.
In this paper, we devote our attention to distributional solutions of \eqref{eq:Intro:EL} belonging to $W^{1,1}(\Omega)$. 
We pursue the question whether such solutions are of locally finite energy which would make the theories of De Giorgi-Nash \cite{EDG, JN58} and Schauder \cite{Schauder1934, Schauder1937, GT} applicable.
In other words, we study functions $u \in W^{1,1}(\Omega)$ for which we have
\[
 \int_{\Omega} Df(Du) \cdot D\varphi \, dx = 0, \quad \forall \varphi \in C_{0}^{\infty}(\Omega).
\]
Notice that if $Df$ is of linear growth, the expression above is well-defined.
Our study is focused on existence and non-existence of distributional solutions $u \in W^{1,1}(\Omega) \setminus W^{1,2}_{\loc}(\Omega)$. Since the energy $\mathcal{E}_f$ of such solutions is infinite, we sometimes call such solutions \emph{non-energetic}. For the Laplace equation it is well-known from Weyl's lemma that any distributional solution to 
\[
 \Delta u = 0 \text{ in } \Omega 
\]
is smooth, proving that no non-energetic solutions exist for the Laplace equation.
The following result due to V. \v{S}ver\'{a}k and X. Yan \cite{VSXY02} treats the vectorial version of \eqref{eq:Intro:EL}: 
\begin{theorem*}[{\cite[Theorem 1 (iii)]{VSXY02}}]
There exists $f \colon \R^{3 \times 3} \to \R$ smooth, uniformly convex such that $|D^2 f| \leq c$ in $\R^{3 \times 3}$ and the equation
\begin{equation*}
 \divergence (Df(Du)) = 0
\end{equation*}
admits a solution in $W^{1,p} \setminus W^{1,2}$ ($1 < p < 2$). This in particular proves non-uniqueness of solutions in $W^{1,p}$.
\end{theorem*}

{The existence and non-existence of non-energetic solutions has already been extensively studied for the $p$-Laplace equation as well as the linear equation
\begin{equation}\label{eq:IntroLinearEq}
 \divergence (A(x) \cdot Du) = 0 \text{ in } \Omega
\end{equation}
where $\lambda I \leq A(x) \leq \Lambda I$ and A is measurable. We give a swift outline of the situation for \eqref{eq:IntroLinearEq}. For the standard regularity theory we refer to the fundamental works by E. De Giorgi and J. Nash \cite{EDG, JN58} (see also \cite{AFVAs, XFRXRO}). Subsequently, J. Serrin \cite{JS64} provided an example of a distributional solution to \eqref{eq:IntroLinearEq} not belonging to $L^{\infty}_{\loc}$.
Indeed, the construction by J. Serrin shows that for every $p \in (1,2)$, there exists a matrix field $A$ such that $\lambda I \leq A(x) \leq \Lambda I$ for some $\lambda$ and $\Lambda$ and with $u \in W^{1,p} \setminus L^{\infty}_{\loc}$ solution to \eqref{eq:IntroLinearEq}. This shows that $u \in W^{1,2}$ is a crucial assumption in the theory by De Giorgi-Nash.
Following this striking achievement, J. Serrin \cite{JS64} conjectured that if $A$ is H\"older continuous, then any $W^{1,1}$ solution to \eqref{eq:IntroLinearEq} belongs to $W^{1,2}_{loc}$. The conjecture was confirmed by R. Hager and J. Ross \cite{RHJR72} for $W^{1,p}$ solutions with $1 < p < 2$ and by H. Brezis \cite{HB08, AA09} for $W^{1,1}$ solutions. The results by H. Brezis are more general than what the conjecture requires, see \cite{HB08, AA09}. For further results regarding very weak solutions to \eqref{eq:IntroLinearEq}, we refer to \cite{TJVMJVAS09}. We also give a brief description of the situation for the $p$-Laplace equation. For the standard regularity theory we refer to \cite{NNU, KUh, LCEvansNewProof, JLL, PTolk}.
In \cite{TICS}, T. Iwaniec and C. Sbordone conjectured that any $W^{1,r}$ solution with $r > \max \{1, p-1 \}$ to the $p$-Laplace equation is $W^{1,p}_{\loc}$ and hence classical in the interior of the domain. The conjecture was confirmed in the case when $p - \delta < r < p$ for some $\delta > 0$ depending on $n$ and $p$ by T. Iwaniec and C. Sbordone \cite{TICS}. 
Recently, in \cite{MCRT}, the conjecture was proved to be false in full generality.}
\\ \\
The main goal of this paper is to establish an optimal condition (see \eqref{eq:AsymptoticConditionLaplacian} below) under which no non-energetic solutions exist. 
Our first main result is the following, which treats more general problems than \eqref{eq:Intro:EL}. It states that if the Euler-Lagrange equation resembles that of a convex quadratic form, then a Weyl-type result holds.
\begin{maintheorem}\label{thm:ExtWeyl}
 Let $n \geq 1$ and $\Omega \subset \R^n$ be a Lipschitz domain. Let $f \in C^1(\R^n)$ be a function so that there exist $0 \leq \gamma < 1$, $A \in \R^{n \times n}$ an elliptic matrix with ellipticity constants $m$ and $M$ as well as a constant $K > 0$ for which we have
 \begin{equation}\label{eq:AsymptoticConditionLaplacian}
 \tag{N}
 |Df(x) - Ax| \leq K(1 + |x|^{\gamma}), \quad \forall x \in \R^n. 
 \end{equation}
Let $h \in L^{p}(\Omega)$ for some $1 < p < n$ and $u \in W^{1,1}(\Omega)$ be a distributional solution of $\divergence (Df(Du)) = h$. Then $u \in W^{1,p^{\star}}_{\loc}(\Omega)$ (where $p^{\star} \coloneqq \frac{np}{n-p}$) with 
 \begin{equation}
  \| u \|_{W^{1, p^{\star}}(\Omega^{\prime})} \leq C \left( 1 + \| u \|_{W^{1, 1}(\Omega)} + \| h \|_{L^{p}(\Omega)} \right) \quad \text{for all $\Omega^{\prime} \Subset \Omega$}
 \end{equation}
 for some constant $C$ depending only on $\Omega$, $\Omega^{\prime}$, $n$, $p$, $\gamma$, $m$, $M$ and $K$.
\end{maintheorem}
\noindent Note that, although we require $A$ to be elliptic, we do not require $f$ to be convex.
Theorem~\ref{thm:ExtWeyl} implies that if $u \in W^{1,1}(\Omega)$ is a distributional solution of \eqref{eq:Intro:EL} and $f$ satisfies \eqref{eq:AsymptoticConditionLaplacian}, then $u \in W^{1,p}_{\loc}(\Omega)$ for all $1 < p < \infty$. \\ \\
Our second main result demonstrates that if \eqref{eq:AsymptoticConditionLaplacian} is not satisfied, then, one can in general not improve the regularity from $W^{1,1}$ to $W^{1,2}_{\loc}$ for solutions to \eqref{eq:Intro:EL}.
\begin{maintheorem}\label{thm:MainFunctionals}
Let $\Omega \subset \R^2$ be a Lipschitz domain and $0 < \lambda < \Lambda < \infty$ be arbitrary. Then there is a smooth and uniformly convex function $f \in C^{\infty}(\R^2)$ with $\lambda I \leq D^2f(x) \leq \Lambda I$, $\forall x \in \R^2$ such that for any $p < p_{\lambda, \Lambda} \coloneqq \frac{2 \sqrt{\Lambda}}{\sqrt{\lambda} + \sqrt{\Lambda}}$ and any $\alpha \in (0,1)$, \eqref{eq:Intro:EL} admits infinitely many solutions $u \in W^{1,p}(\Omega) \cap C^{\alpha}(\overline{\Omega})$ with identical affine boundary datum and 
\begin{equation}\label{eq:NowhereLocallyInTheSpaceW12}
 \int_B |Du|^2 \, dx = + \infty \quad \text{ for any ball } B \subset \Omega.
\end{equation}
\end{maintheorem}
We observe that $p_{\lambda, \Lambda} \to 2$ as $\Lambda \to \infty$ for fixed $\lambda$. 
Moreover, we observe that the exponent $p_{\lambda, \Lambda}$ coincides with the lower critical exponent for 2D linear elliptic equations, see the works by F. Leonetti and V. Nesi \cite{FLVN97}, K. Astala \cite{KA94} and S. Petermichl and A. Volberg \cite{SPAV02} (see also \cite{KATIES01, ODAV03} and \cite[Theorem 1.1]{AstalaFaracoSzekelyhidiEllipticRegularity}). 

As a byproduct, the proof of Theorem~\ref{thm:MainFunctionals} shows that Theorem~\ref{thm:ExtWeyl} is optimal. Indeed, consider the case $n = 2$ and assume that $f \colon \R^2 \to \R$ satisfies \eqref{eq:AsymptoticConditionLaplacian}. Then
\[
 \lim_{a \to - \infty} \frac{\partial_1 f(a,b)}{a} = \lim_{a \to + \infty} \frac{\partial_1 f(a,b)}{a} \quad \forall b \in \R 
 \quad \text{and} \quad
 \lim_{b \to - \infty} \frac{\partial_2 f(a,b)}{b} = \lim_{b \to + \infty} \frac{\partial_2 f(a,b)}{b} \quad \forall a \in \R.
\]
The functions considered in the proof of Theorem~\ref{thm:MainFunctionals} are such that 
\[ 
\lim_{a \to + \infty} \frac{\partial_1 f(a,b)}{a} = \lambda \quad \text{and} \quad \lim_{a \to - \infty} \frac{\partial_1 f(a,b)}{a} = \Lambda \quad \forall b \in \R
\]
and satisfy \eqref{eq:AsymptoticConditionLaplacian} with $\gamma = 1$ and some elliptic matrix $A$.
This demonstrates that the assumptions of Theorem~\ref{thm:ExtWeyl} cannot be relaxed.

Our proof of Theorem~\ref{thm:MainFunctionals} uses the \emph{convex integration method} invented by S. M\"uller and V. \v{S}ver\'{a}k \cite{SMVS} and ever since used to attack a wide range of problems related to wild solutions of PDE's, compare for instance \cite{LSP, COR, GMTDI, HRT, CJPJRT, ST}. \emph{Staircase laminates}, which were introduced by D. Faraco in \cite{FaracoMiltonConjecture} (and subsequently further developed and used in e.g. \cite{ContiFaracoMaggiCounterexamples, ContiFaracoMaggiMullerRankOne, AstalaFaracoSzekelyhidiEllipticRegularity, XFRRT}) are a central part of the proof. Such techniques were recently used in \cite{MCRT} to build very weak solutions to the $p$-Laplace equation. 
Finally, our last main result reveals that under a \emph{pinching condition} on $f$ and \emph{one-sided bounds} on solutions $u$ to \eqref{eq:Intro:EL}, regularity can be improved from $W^{1,1}$ to $W^{1,2}_{\loc}$.
\begin{maintheorem}\label{thm:OneSidedBound}
 Let $n \geq 1$ and $\Omega \subset \R^n$ be a Lipschitz domain. Let $f \in C^2(\R^n)$ be a convex function so that there exist $0 < \lambda < \Lambda < \infty$ for which we have $\lambda I \leq D^2f(x) \leq \Lambda I$ for all $x \in \R^n$. Assume that the following \emph{pinching condition} 
 \begin{equation}\label{eq:PinchingCondition}
 \Lambda^2(1 - n^{-1}) < \lambda^2
 \end{equation}
 holds. Let $u \in W^{1,1}(\Omega)$ be a distributional solution of $\divergence(Df(Du)) = h$ where $h \in L^1(\Omega)$. Assume that there exist $\sigma_i \in \{ -1 , 1\}$ and $L_i \in \R$ such that
 \begin{equation}\label{eq:OneSidedBoundCondition}
  \sigma_i \partial_i u \geq L_i \quad \text{a.e. in $\Omega$ for all $i = 1, \ldots, n$.} 
 \end{equation}
 Then $u \in W^{1,2}_{\loc}(\Omega)$ with
 \begin{equation}\label{eq:QuantitativeEstimateInOneSidedBoundsTheorem}
 \| u \|_{W^{1,2}(\Omega^{\prime})} \leq C \left( \| Du \|_{L^1(\Omega; \R^n)} + \| h \|_{L^1(\Omega)} \right) \left(\| u \|_{L^1(\Omega)} + \| (L_1, \ldots, L_n) \|_{\ell^2} \right) \quad \text{for all $\Omega^{\prime} \Subset \Omega$}
 \end{equation}
 for some constant $C$ depending only on $\Omega$, $\Omega^{\prime}$, $n$, $\Lambda$ and $\lambda$.
\end{maintheorem}
A comparable result for the $p$-Laplace equation is already available. We refer the reader to \cite[Theorem 1.1]{MCRT}, the proof of which is close to that of Theorem~\ref{thm:OneSidedBound}.
Theorem~\ref{thm:OneSidedBound} leaves open whether the pinching condition \eqref{eq:PinchingCondition} is necessary.

\subsection{Plan of the paper}
In Section~\ref{sec:Notation}, we introduce some notation. We prove Theorem~\ref{thm:ExtWeyl} in Section \ref{sec:ProofExtWeyl}. In the next section we prove Theorem~\ref{thm:OneSidedBound}. The remainder of the paper is dedicated to proving Theorem~\ref{thm:MainFunctionals}. We present some preliminaries of convex integration in Section~\ref{sec:Prel-CI}. In Section~\ref{sec:Geometric} we do the necessary geometric construction. In Section~\ref{sec:Scheme}, we give a detailed description of our convex integration scheme. The three last sections prove Theorem~\ref{thm:MainFunctionals}. 

\section*{Acknowledgements}
The author was supported by the SNSF Grant 182565 and by the Swiss State Secretariat for Education, Research and lnnovation (SERI) under contract number MB22.00034.
The author thanks Riccardo Tione and Maria Colombo for suggesting to investigate the topic of very weak solutions to Euler-Lagrange equations. Their suggestions have been of great value. This paper is an outcome of improvements made to results proved in the master thesis of the author carried out at EPFL in fall 2021 \cite{CJMT}. The advice from Xavier Fern\'{a}ndez-Real regarding Theorem~\ref{thm:ExtWeyl} is gratefully acknowledged. Furthermore, the author would like to thank the anonymous referee for careful reading of the manuscript and useful suggestions which contributed to improve the paper.

\section{Notation}\label{sec:Notation}

In this section, we introduce some notation. Throughout this paper, $\Omega \subset \R^n$ denotes any bounded Lipschitz domain in $\R^n$. We will use $| \cdot |$ to denote the $n$-dimensional Lebesgue measure in $\R^n$. We also use this notation to denote the Euclidian norm i.e.
\[
 |x| = |(x_1, \ldots, x_n)| = \left( \sum_{i = 1}^{n} x_i^2 \right)^{\frac{1}{2}}, \quad x \in \R^n.
\]
For any map $u \colon \Omega \to \R^m$ the $j$-th component is denoted by $u^j$.
A continuous Lipschitz map $u \colon \Omega \to \R^m$ will be called piecewise affine if there is a collection $\mathcal{F}$ of open convex mutually disjoint sets covering $\Omega$ up to a set of Lebesgue measure 0 such that for all $U \in \mathcal{F}$, $u |_{U}$ is affine.
Let $n,m \geq 2$. The set of real-valued matrices with $n$ rows and $m$ columns is denoted by $\R^{n \times m}$. 
A matrix $A \in \R^{n \times n}$ is said to be elliptic if there exists constants $0 < m \leq M$ such that $mI \leq A \leq MI$ i.e.
\begin{equation}
m |x|^2 \leq Ax \cdot x \leq M |x|^2 \quad \forall x \in \R^n.
\end{equation}
The constants $m$ and $M$ are called the ellipticity constants of $A$.
The set of real-valued symmetric matrices with $n$ rows and $n$ columns is denoted by $\R^{n \times n}_{\sym}$. 
We say that two matrices $A, B \in \R^{n \times m}$ are rank-one connected if $\rank{A - B} \leq 1$. 
If $A$ and $B$ are rank-one connected, the segment 
\[
 [A,B] \coloneqq \{ tA + (1 - t)B : t \in [0,1] \}
\]
is called a rank-one segment. A function is called rank-one convex if it is convex along any rank-one segment.
We let $\rho \in C^{\infty}_{0}(\R^n)$ be a function such that $\operatorname{supp}(\rho) \subset B_1$ and $\int_{B_1} \rho \, dx = 1$. We define $\rho_{\eps} \in C^{\infty}_{0}(\R^n)$ as
\[
 \rho_{\eps}(x) \coloneqq \frac{1}{\eps^n} \rho \left( \frac{x}{\eps} \right).
\]
For any $u \in W^{1,1}(\Omega)$, we denote its convolution with $\rho_{\eps}$ as $u_{\eps} \coloneqq u \ast \rho_{\eps}$. Note that $u_{\eps}$ is only well-defined in $\{ x \in \Omega: \dist{x}{\Omega^c} > \eps\}$. To overcome this restriction in the case where $u$ has affine boundary datum, we extend $u$ in an affine manner to $\R^n$. Then $u_{\eps}$ is well-defined on the whole $\Omega$. This is also true for $(Du)_{\eps}$. Finally, given $x_1, \ldots, x_N$, we define
\[
 \sum_{k = i}^{j} x_k \coloneqq 0 \quad \text{and} \quad \prod_{k = i}^{j} x_k \coloneqq 1 \quad \text{whenever } j < i.
\]

\section{Proof of Theorem~\ref{thm:ExtWeyl}}\label{sec:ProofExtWeyl}
In this section, we prove Theorem~\ref{thm:ExtWeyl}. We first need a preliminary lemma. We recall that for any $p \in [1,n)$, we define $p^{\star}$ as $p^{\star} = \frac{np}{n - p}$.
\begin{lemma}\label{lemma:PrelLemmaElliptic}
 Let $1 < p < n, 1 < q < \infty$ and $A \in \R^{n \times n}$ be an elliptic matrix with ellipticity constants $m$ and $M$. If $u \in W^{1,1}(\Omega)$ is a distributional solution of
 \begin{equation*}
  \divergence (A Du) = h + \divergence G
 \end{equation*}
 for some $h \in L^p(\Omega)$ and $G \in L^q(\Omega; \R^n)$. Then $u \in W^{1, \min(p^{\star}, q)}_{\loc}(\Omega)$ with 
 \begin{equation}\label{eq:ExtWeylEstimateInStatement}
  \| u \|_{W^{1, \min(p^{\star}, q)}(\Omega^{\prime})} \leq C ( \| u \|_{W^{1, 1}(\Omega)} + \| h \|_{L^{p}(\Omega)} + \| G \|_{L^q(\Omega; \R^n)} ) \quad \text{for all $\Omega^{\prime} \Subset \Omega$}
 \end{equation}
 for some constant $C$ depending only on $\Omega$, $\Omega^{\prime}$, $n$, $p$, $q$, $m$ and $M$.
\end{lemma}

\begin{proof}
 \textbf{Step 1: }Assume that $A$ is the identity matrix. Then $u$ solves $\Delta u = h + \divergence G$. Let $B \subset \Omega$ be an arbitrary ball. By \cite[Theorem 9.9]{GT} and \cite[Section 2 and Eq. (2.9)]{TICS}, there exists $u_B \in W^{1, \min(p^{\star}, q)}(B)$ solving $\Delta u_B = h + \divergence G$ on $B$ with 
\begin{equation}\label{eq:WeylExtLemmaHodgePlusCZ}
 \| u_B \|_{W^{1, \min(p^{\star}, q)}(B)} \leq c (\| h \|_{L^p(B)} + \| G \|_{L^q(B; \R^n)})
\end{equation}
 for some $c \geq 1$ depending only on $B$, $n$, $p$ and $q$.
 Hence $u - u_B$ is harmonic on $B$. Thus $u \in W^{1, \min(p^{\star}, q)}_{\loc}(B)$. 
 By Cald\'eron-Zygmund estimates, for any $B^{\prime} \Subset B \subset \Omega$
 \begin{equation}\label{eq:WeylExtLemmaOnlyCZ}
  \| u - u_B \|_{W^{1, \min(p^{\star}, q)}(B^{\prime})} \leq c^{\prime} \| u - u_B \|_{W^{1, 1}(B)}.
 \end{equation}
 for some $c^{\prime}$ depending only on $B$, $B^{\prime}$, $n$, $p$ and $q$.
 Therefore,
  \begin{align}
 \begin{split}\label{eq:EstimateInaCompactlyEmbeddedBall}
  \| u \|_{W^{1, \min(p^{\star}, q)}(B^{\prime})} &\leq \| u - u_B \|_{W^{1, \min(p^{\star}, q)}(B^{\prime})} + \| u_B \|_{W^{1, \min(p^{\star}, q)}(B^{\prime})} \\
  \overset{\eqref{eq:WeylExtLemmaOnlyCZ}}&{\leq} c^{\prime} \| u - u_B \|_{W^{1, 1}(B)} + \| u_B \|_{W^{1, \min(p^{\star}, q)}(B^{\prime})} \\
  &\leq 2 c^{\prime} ( \| u \|_{W^{1, 1}(B)} + \| u_B \|_{W^{1, \min(p^{\star}, q)}(B)} ) \\
  \overset{\eqref{eq:WeylExtLemmaHodgePlusCZ}}&{\leq} 2 c c^{\prime} (\| u \|_{W^{1, 1}(B)} + \| h \|_{L^p(B)} + \| G \|_{L^q(B; \R^n)} ).
  \end{split}
 \end{align}
 Since $u \in W^{1, \min(p^{\star}, q)}_{\loc}(B)$ for any ball $B \subset \Omega$, we conclude that $u \in W^{1, \min(p^{\star}, q)}_{\loc}(\Omega)$. By virtue of \eqref{eq:EstimateInaCompactlyEmbeddedBall}, for any $\Omega^{\prime} \Subset \Omega$, 
\begin{equation}\label{eq:PrelLemmaEstimateForLaplacian}
  \| u \|_{W^{1, \min(p^{\star}, q)}(\Omega^{\prime})} \leq C  (\| u \|_{W^{1, 1}(\Omega)} + \| h \|_{L^p(\Omega)} + \| G \|_{L^q(\Omega; \R^n)} )
\end{equation}
for some constant $C$ depending only on $\Omega$, $\Omega^{\prime}$, $n$, $p$ and $q$.
 \\
 \textbf{Step 2: }Let us now assume that $A$ is an arbitrary symmetric positive definite matrix. Then there exists an invertible matrix $R \in \R^{n \times n}$ such that $A = R^T I R$. 
 Let $R^{-T} \Omega \coloneqq \{ R^{-T} x : x \in \Omega \}$ and define $v \in W^{1,1}(R^{-T} \Omega)$ as
 \begin{equation*}
  v(x) \coloneqq u(R^T x).
 \end{equation*}
 Then $v$ solves $\Delta v = \tilde{h} + \divergence \widetilde{G}$ where $\tilde{h} \in L^p(R^{-T} \Omega)$ and $\widetilde{G} \in L^q(R^{-T} \Omega; \R^n)$ are defined as
 \begin{equation*}
  \tilde{h}(x) \coloneqq h(R^{T} x) \quad \text{and} \quad \widetilde{G}(x) \coloneqq R^{-T} G(R^{T} x).
 \end{equation*}
 By Step 1, $v \in W^{1, \min(p^{\star}, q)}_{\loc}(R^{-T} \Omega)$. We conclude that $u \in W^{1, \min(p^{\star}, q)}_{\loc}(\Omega)$. By \eqref{eq:PrelLemmaEstimateForLaplacian}, for any $\Omega^{\prime} \Subset \Omega$, 
 \begin{equation}\label{eq:PrelLemmaEstimateForSymetricMatrices}
  \| u \|_{W^{1, \min(p^{\star}, q)}(\Omega^{\prime})} \leq C  (\| u \|_{W^{1, 1}(\Omega)} + \| h \|_{L^p(\Omega)} + \| G \|_{L^q(\Omega; \R^n)} )
 \end{equation}
 for some constant $C$ depending only on $\Omega$, $\Omega^{\prime}$, $n$, $p$, $q$, $m$ and $M$.
 \\
 \textbf{Step 3: }Assume that $A$ is only known to be elliptic. Let $A_{\sym}$ be the symmetric part of $A$ is defined as $A_{\sym} = \frac{1}{2}(A + A^T)$. Observe that
 \begin{equation*}
  \divergence (A_{\sym} Du) = \divergence (A Du) = h + \divergence G.
 \end{equation*}
 By Step 2, $u \in W^{1, \min(p^{\star}, q)}_{\loc}(\Omega)$ and \eqref{eq:PrelLemmaEstimateForSymetricMatrices} holds. This ends the proof of Lemma~\ref{lemma:PrelLemmaElliptic}.
\end{proof}

We can now prove Theorem~\ref{thm:ExtWeyl}.
\begin{proof}[Proof of Theorem~\ref{thm:ExtWeyl}]
We begin with the case when $0 < \gamma < 1$. We find that
\begin{equation*}
\divergence(A Du) = \divergence (\underbrace{A Du - Df(Du)}_{=:  G }) + \divergence (Df(Du)) = \divergence  G  + h.
\end{equation*}
Since 
\begin{equation}\label{eq:BoundsOnG}
| G | \leq |A Du - Df(Du)| \leq K(1 + |Du|^{\gamma}),
\end{equation}
we find that $ G  \in L^{\sfrac{1}{\gamma}}(\Omega; \R^n)$ with
\begin{equation}
 \| G \|_{L^{\sfrac{1}{\gamma}}(\Omega; \R^n)} \leq c_{K,\gamma} (1 + \| Du \|_{L^1(\Omega; \R^n)})
\end{equation}
for some constant $c_{K,\gamma}$ depending only on $K$ and $\gamma$.
Therefore by Lemma~\ref{lemma:PrelLemmaElliptic}, $u \in W^{1, {\min(\sfrac{1}{\gamma}, p^{\star})}}_{\loc}(\Omega)$. Due to \eqref{eq:ExtWeylEstimateInStatement}, for any $\Omega^{\prime} \Subset \Omega$, there exists a constant $c$ depending only $\Omega$, $\Omega^{\prime}$, $n$, $p$, $m$, $M$ and $\gamma$ such that
\begin{equation}\label{eq:EstimateForTheFirstIterationInTheoremA}
\begin{split}
 \| u \|_{W^{1, {\min(\sfrac{1}{\gamma}, p^{\star})}}(\Omega^{\prime})} &\leq c (\| u \|_{W^{1, 1}(\Omega)} + \| h \|_{L^p(\Omega)} + \| G \|_{L^{\sfrac{1}{\gamma}}(\Omega; \R^n)} ) \\
 &\leq c (\| u \|_{W^{1, 1}(\Omega)} + \| h \|_{L^p(\Omega)} + c_{K,\gamma} (1 + \| Du \|_{L^1(\Omega; \R^n)}) ) \\
 &\leq 2 c c_{K,\gamma} (1 + \| u \|_{W^{1, 1}(\Omega)} + \| h \|_{L^p(\Omega)}).
\end{split}
\end{equation}
Thus, by \eqref{eq:BoundsOnG}, $G \in L^{\sfrac{1}{\gamma}\min(\sfrac{1}{\gamma}, p^{\star})}_{\loc}(\Omega; \R^n)$. By applying the same argument as above, $u \in W^{1, \min(\sfrac{1}{\gamma^2}, p^{\star})}_{\loc}(\Omega)$ with
\[
 \| u \|_{W^{1, \min(\sfrac{1}{\gamma^2}, p^{\star})}(\Omega^{\prime})} \leq c^{\prime} (1 + \| u \|_{W^{1, 1}(\Omega)} + \| h \|_{L^p(\Omega)}) \quad \text{ for all $\Omega^{\prime} \Subset \Omega$}
\]
for some constant $c^{\prime}$ depending only on $\Omega$, $\Omega^{\prime}$, $n$, $p$, $\gamma$, $m$, $M$ and $K$.
By a bootstrapping argument on a finite increasing sequence of compacts between $\Omega^{\prime}$ and $\Omega$, $u \in W^{1, p^{\star}}_{\loc}(\Omega)$ with  
\begin{equation}\label{eq:FinalQuantitativeEstimateInProofOfExtWeyl}
 \| u \|_{W^{1, p^{\star}}(\Omega^{\prime})} \leq C (1 + \| u \|_{W^{1, 1}(\Omega)} + \| h \|_{L^p(\Omega)})  \quad \text{ for all $\Omega^{\prime} \Subset \Omega$}
\end{equation}
for some constant $C$ depending only on $\Omega$, $\Omega^{\prime}$, $n$, $p$, $\gamma$, $m$, $M$ and $K$.
If $\gamma = 0$, then the argument above gives that $G \in L^{\infty}(\Omega; \R^n)$. Thus $u \in W^{1, p^\star}_{\loc}(\Omega; \R^n)$ and \eqref{eq:FinalQuantitativeEstimateInProofOfExtWeyl} holds. This finishes the proof.
\end{proof}

\section{Proof of Theorem~\ref{thm:OneSidedBound}}\label{sec:ProofOneSided}
In this section, we shall prove Theorem~\ref{thm:OneSidedBound}. First, however, we need a preliminary lemma. We recall the notation $w_{\eps} \coloneqq w \ast \rho_{\eps}$ from Section~\ref{sec:Notation}.

\begin{lemma}\label{lemma:ResultAboutFuncWithOneSidedBounds}
 Let $n \geq 1$, $\Omega \subset \R^n$ and $u \in W^{1,1}(\Omega)$. Assume that for all $i = 1, \ldots, n$, there exist $\sigma_i \in \{ -1, 1\}$ and $L_i \in \R$ such that
 \[
  \sigma_i \partial_i u \geq L_i \text{ a.e. on } \Omega.
 \]
 Then for any $\Omega^{\prime} \Subset \Omega$ and $\eps < \frac{1}{2}\dist{\overline{\Omega}^{\prime}}{\Omega^c}$
 \begin{equation}\label{eq:CommuteEstimatorOfAbsAndConv}
  (|Du|)_{\eps} \leq n^{\sfrac{1}{2}} |Du_{\eps}| + 2 n^{\sfrac{1}{2}} \| (L_1, \ldots, L_n) \|_{\ell^2} \text{ a.e. on } \Omega^{\prime}.
 \end{equation}
\end{lemma}

\begin{proof}
 \textbf{Step 1: } Fix some arbitrary $\Omega^{\prime} \Subset \Omega$ and $\eps < \frac{1}{2}\dist{\overline{\Omega}^{\prime}}{\Omega^c}$. We start by assuming that all derivatives $\partial_i u$ preserve their sign over $\Omega$. By this, we mean that there exists $\sigma_i \in \{ -1, 1\}$ for all $i= 1, \ldots, n$ so that $\partial_i u = \sigma_i |\partial_i u|$. Then, for all $x \in \Omega^{\prime}$
 \begin{align*}
  (|Du|)_{\eps}(x) &= \int_{\Omega} \left( \sum_{i = 1}^{n} |\partial_i u|^2 \right)^{\sfrac{1}{2}} (x-y) \rho_{\eps}(y) \, dy \\
  &\leq \int_{\Omega} \sum_{i = 1}^{n} |\partial_i u (x-y)| \rho_{\eps}(y) \, dy \\
  &= \sum_{i = 1}^n \left| \int_{\Omega} \partial_i u(x-y) \rho_{\eps}(y) \, dy \right| \\
  &\leq n^{\sfrac{1}{2}} \left( \sum_{i = 1}^{n} \left| \int_{\Omega} \partial_i u(x-y) \rho_{\eps}(y) \, dy \right|^2 \right)^{\sfrac{1}{2}} = n^{\sfrac{1}{2}} |Du_{\eps}|(x).
  \end{align*}
  \textbf{Step 2: }
  Once again, fix some arbitrary $\Omega^{\prime} \Subset \Omega$ and $\eps < \frac{1}{2}\dist{\overline{\Omega}^{\prime}}{\Omega^c}$.
  Now we assume that there are $\sigma_i \in \{ -1, 1 \}$ and $L_i \in \R$ such that $\sigma_i \partial_i u \geq L_i$ a.e. on $\Omega$. Define a new function as $\tilde{u}(x) \coloneqq u(x) - \sum_{i = 1}^{n} \sigma_i L_i x_i$. The function $\tilde{u}$ satisfies the conditions of Step 1 and therefore $(|D\tilde{u}|)_{\eps} \leq n^{\sfrac{1}{2}} |D\tilde{u}_{\eps}|$ a.e. on $\Omega^{\prime}$. We conclude that
  \begin{align*}
   (|Du|)_{\eps}(x) &\leq (|D\tilde{u}|)_{\eps}(x) + \| (L_1, \ldots, L_n) \|_{\ell^2} \\
   &\leq n^{\sfrac{1}{2}} |D\tilde{u}_{\eps}|(x) + \| (L_1, \ldots, L_n) \|_{\ell^2} \\
   &\leq n^{\sfrac{1}{2}} |Du_{\eps}|(x) + 2 n^{\sfrac{1}{2}}\| (L_1, \ldots, L_n) \|_{\ell^2}
  \end{align*}
  for all $x \in \Omega^{\prime}$.
\end{proof}

We can now prove Theorem~\ref{thm:OneSidedBound}. The proof is close to the proof of \cite[Theorem 1.1]{MCRT}.
\begin{proof}[Proof of Theorem~\ref{thm:OneSidedBound}]
 \textbf{Step 1: } In this first step, we assume that $n(\Lambda^2 - 2\lambda + 1) < 1$. The reason for this assumption is not clear at this point but will become clear in the proof. Moreover, without loss of generality, we can assume that $\lambda < \Lambda$. Indeed, if $\lambda = \Lambda$, then there is a $\lambda^{\prime} < \Lambda$ such that $n(\Lambda^2 - 2\lambda^{\prime} + 1) < 1$ and $\lambda^{\prime} I \leq D^2 f(x) \leq \Lambda I$ still holds. Notice that, since we assume $\lambda < \Lambda$, we have $\Lambda^2 - 2 \lambda + 1 > \Lambda^2 - 2 \Lambda + 1 = (\Lambda - 1)^2 \geq 0$. In addition, we may, without loss of generality, assume that $Df(0) = 0$.
 We have 
 \[
  \int_{\Omega} Df(Du) \cdot D \varphi \, dx = \int_{\Omega} h \varphi \, dx \quad \forall \varphi \in C^{\infty}_0(\Omega).
 \]
 Let $\Omega^{\prime} \Subset \Omega^{\prime \prime} \Subset \Omega^{\prime \prime \prime} \Subset \Omega$ with $\dist{\overline{\Omega}^{\prime}}{\Omega^c} \leq 2 \dist{\overline{\Omega}^{\prime \prime \prime}}{\Omega^c}$ and let $\eps < \dist{\overline{\Omega}^{\prime \prime \prime}}{\Omega^c}$. We find that 
 \[
  \int_{\Omega^{\prime \prime}} (Df(Du))_{\eps} \cdot D \varphi \, dx = \int_{\Omega^{\prime \prime}} h_{\eps} \varphi \, dx \quad \forall \varphi \in C^{\infty}_0(\Omega^{\prime \prime}).
 \]
Let $\Psi \in C^{\infty}_0(\Omega^{\prime \prime})$ such that $\Psi = 1$ on $\Omega^{\prime}$ and $\| D \Psi \|_{L^{\infty}(\Omega^{\prime \prime}; \R^n)} \leq 5 \dist{\overline{\Omega}^{\prime}}{{\Omega^{\prime \prime}}^c}^{-1}$. Then
\[
  \int_{\Omega^{\prime \prime}} (Df(Du))_{\eps} \cdot D(\Psi u_{\eps}) \, dx = \int_{\Omega^{\prime \prime}} h_{\eps} \Psi u_{\eps} \, dx.
 \]
 which leads to 
 \begin{equation}\label{eq:ConsenguenceOfEL}
  \int_{\Omega^{\prime \prime}} \left((Df(Du))_{\eps} \cdot Du_{\eps}\right) \Psi \, dx = - \int_{\Omega^{\prime \prime}} \left((Df(Du))_{\eps} \cdot D \Psi \right) u_{\eps} \, dx + \int_{\Omega^{\prime \prime}} h_{\eps} \Psi u_{\eps} \, dx.
 \end{equation}
 The left-hand side in \eqref{eq:ConsenguenceOfEL} equals
 \begin{equation}\label{eq:LeftHandSideInConsenguenceOfEL}
  \int_{\Omega^{\prime \prime}} |Du_{\eps}|^2 \Psi \, dx + \int_{\Omega^{\prime \prime}} ((Df(Du))_{\eps} - Du_{\eps}) \cdot Du_{\eps} \Psi \, dx.
 \end{equation}
Observe that
 \begin{align*}
  |Df(Du) - Du|^2 &= |Df(Du)|^2 - 2 Df(Du) \cdot Du + |Du|^2 \\
  &= \underbrace{|Df(Du)|^2}_{\leq \Lambda^2 |Du|^2} - 2 \underbrace{(Df(Du) - Df(0)) \cdot (Du - 0)}_{\geq \lambda |Du|^2} + |Du|^2 \leq (\Lambda^2 - 2 \lambda + 1) |Du|^2 \\
 \end{align*}
 and hence, by Lemma~\ref{lemma:ResultAboutFuncWithOneSidedBounds}, we have
 \begin{align*}
 |(Df(Du))_{\eps} - Du_{\eps}| &\leq |Df(Du) - Du|_{\eps} \leq (\Lambda^2 - 2 \lambda + 1)^{\sfrac{1}{2}}(|Du|)_{\eps} \\
 \overset{\eqref{eq:CommuteEstimatorOfAbsAndConv}}&{\leq} n^{\sfrac{1}{2}}(\Lambda^2 - 2 \lambda + 1)^{\sfrac{1}{2}}|Du_{\eps}| + 2 n^{\sfrac{1}{2}} (\Lambda^2 - 2 \lambda + 1)^{\sfrac{1}{2}}\| (L_1, \ldots, L_n) \|_{\ell^2}.
 \end{align*}
 Therefore
 \begin{equation}\label{eq:SomeQuantity}
 \begin{split}
  \left| \int_{\Omega^{\prime \prime}} ((Df(Du))_{\eps} - Du_{\eps}) \cdot Du_{\eps} \Psi \, dx \right| &\leq n^{\sfrac{1}{2}}(\Lambda^2 - 2 \lambda + 1)^{\sfrac{1}{2}} \int_{\Omega^{\prime \prime}} |Du_{\eps}|^2 \Psi \, dx \\
  & \qquad + 2 n^{\sfrac{1}{2}} (\Lambda^2 - 2 \lambda + 1)^{\sfrac{1}{2}}\| (L_1, \ldots, L_n) \|_{\ell^2} \int_{\Omega^{\prime \prime}} |Du_{\eps}| \Psi \, dx.
 \end{split}
 \end{equation}
 We observe that due to \cite[Eq. (2.6)]{MCRT} and \eqref{eq:OneSidedBoundCondition} there exists a constant ${C = C(n, \Omega, \dist{\overline{\Omega}^{\prime}}{\Omega^c})}$ such that
 \begin{equation}\label{eq:InLInftyBecauseOfOneSidedBounds}
  \| u_{\eps^{\prime}} \|_{L^{\infty}(\Omega^{\prime \prime})} \leq C(\| u_{\eps^{\prime}} \|_{L^1(\Omega^{\prime \prime \prime})} + \| (L_1, \ldots, L_n) \|_{\ell^2}) \quad \forall \eps^{\prime} < \dist{\overline{\Omega}^{\prime \prime \prime}}{\Omega^c}.
 \end{equation}
 Then, by combining \eqref{eq:ConsenguenceOfEL} with \eqref{eq:LeftHandSideInConsenguenceOfEL} and \eqref{eq:SomeQuantity}, we find
 \begin{align*}
  &\left( 1 - n^{\sfrac{1}{2}}(\Lambda^2 - 2 \lambda + 1)^{\sfrac{1}{2}} \right) \int_{\Omega^{\prime \prime}} |Du_{\eps}|^2 \Psi \, dx \\
  &\leq 2 n^{\sfrac{1}{2}} (\Lambda^2 - 2 \lambda + 1)^{\sfrac{1}{2}}\| (L_1, \ldots, L_n) \|_{\ell^2} \int_{\Omega^{\prime \prime}} |Du_{\eps}| \Psi \, dx \\
  &\qquad+ \left| \int_{\Omega^{\prime \prime}} \left((Df(Du))_{\eps} \cdot D \Psi \right) u_{\eps} \, dx \right| + \left| \int_{\Omega^{\prime \prime}} h_{\eps} \Psi u_{\eps} \, dx \right| \\
  &\leq 2 n^{\sfrac{1}{2}} (\Lambda^2 - 2 \lambda + 1)^{\sfrac{1}{2}}\| (L_1, \ldots, L_n) \|_{\ell^2} \| Du_{\eps} \|_{L^1(\Omega^{\prime \prime}; \R^2)} \\
  &\qquad + \| (Df(Du))_{\eps} \|_{L^{1}(\Omega^{\prime \prime}; \R^n)} \| D \Psi \|_{L^{\infty}(\Omega^{\prime \prime}; \R^n)} \| u_{\eps} \|_{L^{\infty}(\Omega^{\prime \prime})} + \| h_{\eps} \|_{L^{1}(\Omega^{\prime \prime})} \| \Psi \|_{L^{\infty}(\Omega^{\prime \prime})} \| u_{\eps} \|_{L^{\infty}(\Omega^{\prime \prime})} \\
  \overset{\eqref{eq:InLInftyBecauseOfOneSidedBounds}}&{\leq} 2 n^{\sfrac{1}{2}} (\Lambda^2 - 2 \lambda + 1)^{\sfrac{1}{2}}\| (L_1, \ldots, L_n) \|_{\ell^2} \| Du_{\eps} \|_{L^1(\Omega^{\prime \prime}; \R^2)} \\
  &\qquad + C \Lambda \| Du \|_{L^1(\Omega; \R^2)} \| D \Psi \|_{L^{\infty}(\Omega^{\prime \prime}; \R^n)} (\| u_{\eps} \|_{L^1(\Omega^{\prime \prime \prime})} + \| (L_1, \ldots, L_n) \|_{\ell^2} ) \\
  &\qquad + C \| h \|_{L^1(\Omega)} (\| u_{\eps} \|_{L^1(\Omega^{\prime \prime \prime})} + \| (L_1, \ldots, L_n) \|_{\ell^2} ) \\
  &\leq 2 n^{\sfrac{1}{2}} (\Lambda^2 - 2 \lambda + 1)^{\sfrac{1}{2}}\| (L_1, \ldots, L_n) \|_{\ell^2} \| Du \|_{L^1(\Omega; \R^2)} \\
  &\qquad + C ( \Lambda \| Du \|_{L^1(\Omega; \R^2)} \| D \Psi \|_{L^{\infty}(\Omega^{\prime \prime}; \R^n)} + \| h \|_{L^1(\Omega)})(\| u \|_{L^1(\Omega)} + \| (L_1, \ldots, L_n) \|_{\ell^2} ) =: C^{\prime}. \\
 \end{align*}
 By lower semi-continuity,
 \[
  \int_{\Omega^{\prime}} |Du|^2 \, dx \leq \int_{\Omega^{\prime \prime}} |Du|^2 \Psi \, dx \leq \liminf_{\eps \to 0} \int_{\Omega^{\prime \prime}} |Du_{\eps}|^2 \Psi \, dx \leq \dfrac{C^{\prime}}{1 - n^{\sfrac{1}{2}}(\Lambda^2 - 2 \lambda + 1)^{\sfrac{1}{2}}}
< \infty. 
 \]
 Thus, $u \in W^{1,2}(\Omega^{\prime})$. Therefore, we conclude $u \in W^{1,2}_{\loc}(\Omega)$ and \eqref{eq:QuantitativeEstimateInOneSidedBoundsTheorem} holds. \\
 \textbf{Step 2:} In this step, we drop the assumption that $n (\Lambda^2 - 2 \lambda + 1) < 1$. Since $u$ solves the equation $\divergence (Df(Du)) = h$, $u$ also solves the equation
 \begin{equation}\label{eq:NormalizedEquation}
  \divergence \left( \frac{1}{\beta}Df(Du) \right) = \frac{1}{\beta} h
 \end{equation}
 for all $\beta > 0$.
 By direct computations, we find that the assumption that $\Lambda^2(1 - n^{-1}) < \lambda^2$ is sufficient to deduce that there exists a $\beta > 0$ such that
 \[
  n \left( \frac{\Lambda^2}{\beta^2} - 2 \frac{\lambda}{\beta} + 1 \right) < 1.
 \]
 We can then apply the arguments from Step 1 to equation \eqref{eq:NormalizedEquation} to conclude the proof. Since the choice of $\beta$ depends only on $n$, $\Lambda$ and $\lambda$, \eqref{eq:QuantitativeEstimateInOneSidedBoundsTheorem} still holds.
\end{proof}

\section{Preliminaires of classical Convex Integration}\label{sec:Prel-CI}
In this section, we introduce some results from classical convex integration, see \cite{SMVS, AstalaFaracoSzekelyhidiEllipticRegularity}.
We introduce the building blocks of the convex integration scheme invented by M\"uller and \v{S}ver\'{a}k and some extensions due to Astala, Faraco and Sz\'ekelyhidi. The content of this subsection comes from \cite{SMVS, AstalaFaracoSzekelyhidiEllipticRegularity}.

\begin{definition}
The barycenter of a probability measure $\nu$ on $\R^{n \times m}$ is $\overline{\nu} = \int_{\R^{n \times m}} X \, d\nu(X) \in \R^{n \times m}$.
\end{definition}

\begin{definition}
A probability measure on $\R^{n \times m}$ $\nu$ is a laminate if $f(\overline{\nu}) \leq \int_{\R^{n \times m}} f(X) \, d\nu(X)$ for all $f \colon \R^{n \times m} \to \R$ rank-one convex.
\end{definition}

\begin{definition}
Assume $\nu = \sum_{j = 1}^r \lambda_{j} \delta_{A_j}$ with $A_j \neq A_k$ whenever $j \neq k$. The measure $\nu^{\prime}$ is said to be obtained from $\nu$ by an elementary splitting if for some $j$ and some $\lambda \in [0,1]$, there exists a rank-one segment $[B_1, B_2]$ containing $A_j$ with $(1 - s)B_1 + sB_2 = A_j$ such that
\[
 \nu^{\prime} = \nu + \lambda \lambda_j \left( (1 - s)\delta_{B_1} + s\delta_{B_2} - \delta_{A_j} \right).
\]
\end{definition}
It follows from the definition of laminates that any probability measure that is obtained by an elementary splitting of a laminate is also a laminate. This justifies the notion of laminates of finite order:

\begin{definition}
 We say that a laminate $\nu$ is a laminate of finite order if there exist probability measures $\nu_1, \ldots, \nu_m$ such that
 $\nu_1 = \delta_A$ for some $A$, $\nu_m = \nu$ and for each $j = 1, \ldots, m-1$, $\nu_{j+1}$ can be obtained from $\nu_j$ by an elementary splitting.
\end{definition}

\begin{lemma}\label{lemma:ApproxLaminatesFiniteOrder}
Let $\nu = \sum_{j = 1}^r \lambda_j \delta_{A_j}$ (with $A_j \in \R^{n \times m}$ for all $j$ and $A_j \neq A_k$ whenever $j \neq k$) be a laminate of finite order and let $A = \overline{\nu}$ be the barycenter of $\nu$. Then for any $\alpha \in (0,1)$, any $\eps > 0$, any $b \in \R^m$ and any $0 < \delta < \frac{1}{2} \min \{ |A_i - A_j| : i \neq j \}$ there is a piecewise affine mapping $u \colon \Omega \to \R^m$ such that
\begin{enumerate}
 \item $\| u - (Ax + b) \|_{C^{\alpha}(\overline{\Omega}; \R^m)} < \eps$;
 \item $u = Ax + b$ on $\partial \Omega$;
 \item $|\{ x \in \Omega : |Du(x) - A_j| < \delta \}| = \lambda_j |\Omega|$ for all $j = 1, \ldots, r$.
\end{enumerate}
If $m = n$ and all the matrices $A_j \in \R^{n \times n}_{\sym}$, then the piecewise affine mapping $u$ can be selected so that there exists $g \in W^{2, \infty}(\Omega)$ for which we have $u = Dg$.
\end{lemma}

\begin{remark}\label{rmk:AboutHowBuildingBlocksAreUsed}
The last sentence in the previous lemma is particularly important in our construction. Indeed, we will use the following which is a consequence of the lemma above: let $\nu = \sum_{j = 1}^r \lambda_j \delta_{A_j}$ be a laminate of finite order with $A_j \in \R^{n \times n}_{\sym}$ for all $j = 1, \ldots, r$, $A_i \neq A_j$ whenever $i \neq j$ and $A = \overline{\nu}$. In addition, for all $j = 1, \ldots, r$, let $U_j \subset \R^{n \times n}_{\sym}$ be an open neighbourhood of $A_j$ in $\R^{n \times n}_{\sym}$. Then for any $\alpha \in (0,1)$, any $\eps > 0$ and any $b \in \R^n$, there exists a piecewise affine mapping $u \colon \Omega \to \R^n$ such that
\begin{enumerate}
 \item $\| u - (Ax + b) \|_{C^{\alpha}(\overline{\Omega}; \R^n)} < \eps$;
 \item $u = Ax + b$ on $\partial \Omega$;
 \item $|\{ x \in \Omega : Du(x) \in U_j \}| = \lambda_j |\Omega|$ for all $j = 1, \ldots, r$.
\end{enumerate}
\end{remark}

\section{Geometric constructions (Part 1 of the proof of Theorem~\ref{thm:MainFunctionals})}\label{sec:Geometric}
In this section, we do the necessary geometric constructions to be able to perform a convex integration scheme in the next section. In Subsection~\ref{subsec:EulerLagrangeAsDifferentialInclusion}, we show that \eqref{eq:Intro:EL} is equivalent to a differential inclusion $Dw \in K_f$. In this first subsection, we also select $f$ depending on the values of $\lambda$ and $\Lambda$. The set $K_f$ is central in the geometric constructions. In Subsection~\ref{subsec:ConstructionAndParametrisationOfLaminates}, we construct laminates of finite order partly supported on $K_f$. Although partly supported in the set $K_f$, it is not clear how to use these laminates together with Lemma~\ref{lemma:ApproxLaminatesFiniteOrder} (see also Remark~\ref{rmk:AboutHowBuildingBlocksAreUsed}) in order to construct a functioning convex integration scheme. Therefore, in Subsection~\ref{subsec:InterpolationMaps}, we define interpolation maps. In Subsection~\ref{subsec:ConstructionParametrisationInterpolationLaminates}, we use these maps to modify the already constructed laminates so that they can be used to construct a functioning convex integration scheme.

\subsection{Equation \eqref{eq:Intro:EL} as a differential inclusion and choice of $f$}\label{subsec:EulerLagrangeAsDifferentialInclusion}
The aim of this subsection is to show that \eqref{eq:Intro:EL} is equivalent to a differential inclusion. Define
\[
 K_f \coloneqq 
 \left\{
 \begin{pmatrix}
 a & b \\
 \partial_2 f(a,b) & - \partial_1 f(a,b) \\
 \end{pmatrix}
 :
 a, b \in \R
 \right\}.
\]
Assume $u \in W^{1,1}(\Omega)$ solves \eqref{eq:Intro:EL}. Define
\[
 J \coloneqq
 \begin{pmatrix}
 0 & 1 \\
 -1 & 0 \\
 \end{pmatrix}.
\]
Since $u$ solves \eqref{eq:Intro:EL}, $\rotational{JDf(Du)} = 0$. Due to this, there is $\tilde{u} \in W^{1,1}(\Omega)$ such that $D \tilde{u} = JDf(Du)$. Define $w \coloneqq (u, \tilde{u}) \in W^{1,1}(\Omega; \R^2)$. Then $Dw(x) \in K_f$ for a.e. $x\in \Omega$. Now instead assume that for some $w \in W^{1,1}(\Omega; \R^2)$, we have $Dw(x) \in K_f$ for a.e. $x \in \Omega$. Then $\divergence(Df(Dw^1)) = \rotational{Dw^2} = 0$. Therefore, equation \eqref{eq:Intro:EL} is equivalent to the differential inclusion $Dw(x) \in K_f$ for a.e. $x \in \Omega$. Given $0 < \lambda < \Lambda$ as in Theorem~\ref{thm:MainFunctionals}, let $\varphi \in C^{\infty}(\R)$ be any function such that $0 < \lambda < \inf_{a \in \R} \varphi^{\prime \prime}(a) \leq \sup_{a \in \R} \varphi^{\prime \prime}(a) < \Lambda$,
\[
 \lim_{a \to - \infty} \frac{\varphi^{\prime}(a)}{a} = \Lambda \quad \text{and} \quad \lim_{a \to + \infty} \frac{\varphi^{\prime}(a)}{a} = \lambda.
\]
We assume, without loss of generality, that $\lambda < 1 < \Lambda$.
Then, choose $f \colon \R^2 \to \R$ to be $f(a,b) = \varphi(a) + \sfrac{b^2}{2}$. The set $K_f$ becomes
\[
K_f =
 \left\{
 \begin{pmatrix}
 a & b \\
 b & - \varphi^{\prime}(a) \\
 \end{pmatrix}
 :
 a, b \in \R
 \right\}
 \subset
 \R^{2 \times 2}_{\sym}.
\]

\subsection{Construction and parametrisation of laminates of finite order}\label{subsec:ConstructionAndParametrisationOfLaminates}
In this subsection, we build laminates of finite order depending on a set of parameters. The space of parameters is defined as 
\[
 \mathcal{P} \coloneqq (1,2) \times (1,2) \times (-1,1) \subset \R^3.
\]
We define the maps $a_0^+, a_0^-, b \colon \mathcal{P} \to \R$ as follows. Let $P = (a_0^+, a_0^-, b) \in \mathcal{P}$ and define
\[
 a_0^{+}(P) \coloneqq a_0^+, \quad a_0^{-}(P) \coloneqq a_0^-, \quad \text{and} \quad b(P) \coloneqq b.
\]
In other words, $a_0^+$ is the projection from $\mathcal{P}$ into the first coordinate, $a_0^-$ is the projection from $\mathcal{P}$ into the second coordinate and so on. Let $1 < r < 2$ be a constant to be determined later (see Lemma~\ref{lemma:BehaviourOfLambda3}). For each $P \in \mathcal{P}$, we define the sequences $\{ a_k^{+}(P) \}_{k = 1}^{\infty}$ and $\{ a_k^{-}(P) \}_{k = 1}^{\infty}$ as
\[
 a_k^+(P) \coloneqq a_0^+(P) + r^k \quad \text{and} \quad a_k^-(P) \coloneqq a_0^-(P) + \sqrt{\frac{\lambda}{r \Lambda}} r^k
\]
for all $k \geq 1$. The justification for this notation is that $a_k^-(P)$ will always be used with a minus sign in front of it and $a_k^+(P)$ without a minus sign. 
For each $P \in \mathcal{P}$ and each $i \geq 1$, we define
\[
 A_i(P) \coloneqq 
 \begin{pmatrix}
 a_i^+(P) & b(P) \\
 b(P) & - \varphi^{\prime}(- a_i^{-}(P))
 \end{pmatrix}.
\]
We split the matrix $A_i(P)$ along the rank-one direction
\[
 \begin{pmatrix}
 0 & 0 \\
 0 & x
 \end{pmatrix}, \quad x \in \R
\]
with endpoints $B_i(P)$ and $C_i(P)$ where
\[
 B_i(P) \coloneqq 
 \begin{pmatrix}
 a_i^+(P) & b(P) \\
 b(P) & - \varphi^{\prime}(a_i^{+}(P))
 \end{pmatrix}
 \quad 
 \text{and}
 \quad
 C_i(P) \coloneqq 
 \begin{pmatrix}
 a_i^+(P) & b(P) \\
 b(P) & - \varphi^{\prime}(- a_{i+1}^{-}(P))
 \end{pmatrix}.
\]
We have $A_i(P) = \lambda_i^1(P) B_i(P) + (1 - \lambda_i^1(P)) C_i(P)$ with
\[
 \lambda_i^1(P) \coloneqq \frac{\varphi^{\prime}(- a_i^{-}(P)) - \varphi^{\prime}(- a_{i+1}^{-}(P))}{\varphi^{\prime}(a_i^{+}(P)) - \varphi^{\prime}(- a_{i+1}^{-}(P))}.
\]
Then we split the matrix $C_i(P)$ along the rank-one direction
\[
 \begin{pmatrix}
 x^{\prime} & 0 \\
 0 & 0 \\
 \end{pmatrix},
 \quad x^{\prime} \in \R
\]
with endpoints $D_i(P)$ and $E_i(P)$ where
\[
 D_i(P) \coloneqq 
 \begin{pmatrix}
 - a_{i+1}^-(P) & b(P) \\
 b(P) & - \varphi^{\prime}(- a_{i+1}^{-}(P))
 \end{pmatrix}
 \quad 
 \text{and}
 \quad
 E_i(P) \coloneqq 
 \begin{pmatrix}
 a_{i+1}^+(P) & b(P) \\
 b(P) & - \varphi^{\prime}(- a_{i+1}^{-}(P))
 \end{pmatrix}.
\]
We have $C_i(P) = \frac{\lambda_i^2(P)}{1 - \lambda_i^1(P)} D_i(P) + \frac{\lambda_i^3(P)}{1 - \lambda_i^1(P)} E_i(P)$ where 
\[
 \lambda_i^2(P) \coloneqq \frac{\varphi^{\prime}(a_i^{+}(P)) - \varphi^{\prime}(- a_{i}^{-}(P))}{\varphi^{\prime}(a_i^{+}(P)) - \varphi^{\prime}(- a_{i+1}^{-}(P))} \frac{a_{i+1}^{+}(P) - a_{i}^{+}(P)}{a_{i+1}^{-}(P) + a_{i+1}^{+}(P)}
\]
and
\[
 \lambda_i^3(P) \coloneqq \frac{\varphi^{\prime}(a_i^{+}(P)) - \varphi^{\prime}(- a_{i}^{-}(P))}{\varphi^{\prime}(a_i^{+}(P)) - \varphi^{\prime}(- a_{i+1}^{-}(P))} \frac{a_{i+1}^{-}(P) + a_{i}^{+}(P)}{a_{i+1}^{-}(P) + a_{i+1}^{+}(P)}.
\]
Notice that $E_i(P) = A_{i+1}(P)$. By the construction above, we have
\[
 A_i(P) = \lambda_i^1(P) B_i(P) + \lambda_i^2(P) D_i(P) + \lambda_i^3(P) E_i(P) = \lambda_i^1(P) B_i(P) + \lambda_i^2(P) D_i(P) + \lambda_i^3(P) A_{i+1}(P)
\]
and 
\begin{equation}\label{eq:LaminateOfFiniteOrderMu}
 \mu^{(i)}(P) \coloneqq \lambda_i^1(P) \delta_{B_i(P)} + \lambda_i^2(P) \delta_{D_i(P)} + \lambda_i^3(P) \delta_{A_{i+1}(P)}
\end{equation}
is a laminate of finite order with barycenter $A_i(P)$. For an illustration of the matrices $A_i(P)$, $B_i(P)$, $C_i(P)$, $D_i(P)$ and $E_i(P)$, we refer to Figure~\ref{fig:IllustrationsLaminates}. We end this subsection by a result about the behaviour of $\lambda_i^3$ as $i$ is large. We recall that $p$ is a fixed arbitrary quantity in Theorem~\ref{thm:MainFunctionals}, which satisfies 
\[
 p < p_{\lambda, \Lambda} = \frac{2 \sqrt{\Lambda}}{\sqrt{\lambda} + \sqrt{\Lambda}}.
\]
\begin{lemma}\label{lemma:BehaviourOfLambda3}
For any $0 < \delta < \frac{1}{10}(p_{\lambda, \Lambda} - p)$, there exists $1 < r < 2$ and an integer $I_0$ such that for all $i \geq I_0$
\begin{equation}\label{eq:Asymp-3}
\frac{1}{r^{2}} < \inf_{P \in \mathcal{P}} \lambda_i^3(P) \leq \sup_{P \in \mathcal{P}} \lambda_i^3(P) < \frac{1}{r^{p + 2 \delta}}.
\end{equation}
\end{lemma}
\begin{proof}
By standard calculations,
\[
\limsup_{i \to \infty} \sup_{P \in \mathcal{P}} \lambda_i^3(P) \leq \left( \frac{\lambda \sqrt{r} + \sqrt{\Lambda \lambda} }{\lambda \sqrt{r} + \sqrt{\Lambda \lambda} r} \right)^2 \leq \liminf_{i \to \infty} \inf_{P \in \mathcal{P}} \lambda_i^3(P).
\]
We will prove that 
\begin{equation}\label{eq:LogRLimitEquation}
 \lim_{r \to 1} \log_r \left[  \left( \frac{\lambda \sqrt{r} + \sqrt{\Lambda \lambda} r}{\lambda \sqrt{r} + \sqrt{\Lambda \lambda} } \right)^2 \right] = \frac{2 \sqrt{\Lambda}}{\sqrt{\lambda} + \sqrt{\Lambda}} = p_{\lambda, \Lambda}.
\end{equation}
Indeed,
\[
 \log_r \left[  \left( \frac{\lambda \sqrt{r} + \sqrt{\Lambda \lambda} r}{\lambda \sqrt{r} + \sqrt{\Lambda \lambda} } \right)^2 \right] = \dfrac{\log \left[  \left( \frac{\lambda \sqrt{r} + \sqrt{\Lambda \lambda} r}{\lambda \sqrt{r} + \sqrt{\Lambda \lambda} } \right)^2 \right]}{\log r} = \dfrac{2 \log \left( \frac{\lambda \sqrt{r} + \sqrt{\Lambda \lambda} r}{\lambda \sqrt{r} + \sqrt{\Lambda \lambda} } \right) }{\log r} 
\]
so that by Bernoulli-L'H\^opital's rule
\begin{equation}\label{eq:ApplicationBernoulliHopital}
 \lim_{r \to 1} \log_r \left[ \left( \frac{\lambda \sqrt{r} + \sqrt{\Lambda \lambda} r}{\lambda \sqrt{r} + \sqrt{\Lambda \lambda} } \right)^2 \right] = \lim_{r \to 1} 2 \dfrac{ \frac{\lambda \sqrt{r} + \sqrt{\Lambda \lambda}}{\lambda \sqrt{r} + \sqrt{\Lambda \lambda} r }}{r^{-1}} \dfrac{d}{dr}\left( \frac{\lambda \sqrt{r} + \sqrt{\Lambda \lambda} r}{\lambda \sqrt{r} + \sqrt{\Lambda \lambda} } \right).
\end{equation}
By direct computations,
\begin{equation}\label{eq:DerivativeComputationDR}
 \dfrac{d}{dr}\left( \frac{\lambda \sqrt{r} + \sqrt{\Lambda \lambda} r}{\lambda \sqrt{r} + \sqrt{\Lambda \lambda} } \right) = \dfrac{\left( \frac{\lambda \sqrt{\Lambda \lambda} }{2 \sqrt{r}} + \frac{1}{2} \lambda \sqrt{\Lambda \lambda r} + \Lambda \lambda \right)}{(\lambda \sqrt{r} + \sqrt{\Lambda \lambda})^2}.
\end{equation}
Therefore, by combining \eqref{eq:ApplicationBernoulliHopital} and \eqref{eq:DerivativeComputationDR},
\begin{align*}
 \lim_{r \to 1} \log_r \left[ \left( \frac{\lambda \sqrt{r} + \sqrt{\Lambda \lambda} r}{\lambda \sqrt{r} + \sqrt{\Lambda \lambda} } \right)^2 \right] &= 2 \frac{\lambda + \sqrt{\Lambda \lambda}}{\lambda + \sqrt{\Lambda \lambda} }\dfrac{\left( \frac{\lambda \sqrt{\Lambda \lambda} }{2} + \frac{1}{2} \lambda \sqrt{\Lambda \lambda} + \Lambda \lambda \right)}{(\lambda + \sqrt{\Lambda \lambda})^2}\\
 &= 2 \dfrac{\left( \lambda \sqrt{\Lambda \lambda} + \Lambda \lambda \right)}{(\lambda + \sqrt{\Lambda \lambda})^2} = 2 \dfrac{\left( \lambda + \sqrt{\Lambda \lambda} \right) \sqrt{\Lambda \lambda}}{(\lambda + \sqrt{\Lambda \lambda})^2} \\
 &= 2 \dfrac{ \sqrt{\Lambda \lambda}}{(\lambda + \sqrt{\Lambda \lambda})} = 2 \dfrac{ \sqrt{\Lambda \lambda}}{\sqrt{\lambda} (\sqrt{\lambda} + \sqrt{\Lambda})} = \dfrac{2 \sqrt{\Lambda}}{\sqrt{\lambda} + \sqrt{\Lambda}} = p_{\lambda, \Lambda}.
\end{align*}
This proves \eqref{eq:LogRLimitEquation}.
Thus, there exists an $1 < r < 2$ such that 
\[
 2 > \log_r \left[  \left( \frac{\lambda \sqrt{r} + \sqrt{\Lambda \lambda} r}{\lambda \sqrt{r} + \sqrt{\Lambda \lambda} } \right)^2 \right] > p + 2 \delta 
\]
which implies
\[
 r^2 > \left( \frac{\lambda \sqrt{r} + \sqrt{\Lambda \lambda} r}{\lambda \sqrt{r} + \sqrt{\Lambda \lambda} } \right)^2 > r^{p + 2 \delta}.
\]
We conclude \eqref{eq:Asymp-3}.
\end{proof}
From now on, we let $\delta$ be a fixed quantity satisfying $0 < \delta < \frac{1}{10}(p_{\lambda, \Lambda} - p)$. The quantities $r$ and $I_0$ denote those given by the above lemma. In addition, we let $C_{r,p} \geq 0$ denote a constant for which we have
\begin{equation}\label{eq:DefinitionOfCRP}
 j r^{-p(j - 1)} \leq C_{r,p} \quad \forall j \geq 1.
\end{equation}

\subsection{Definition and properties of interpolation maps}\label{subsec:InterpolationMaps}
The goal of this subsection is to define interpolation maps and prove useful properties of these maps. For any $i \geq 1$, $t \in (0,1)$, we define $\Phi_{i,t}^1, \Phi_{i,t}^2 \colon \mathcal{P} \to \R^{2 \times 2}_{\operatorname{sym}}$ as
\[
 \Phi_{i,t}^1(P) = tB_i(P) + (1 - t)A_i(P) \quad \text{and} \quad \Phi_{i,t}^2(P) = tD_i(P) + (1 - t)A_{i+1}(P).
\]
\begin{figure}
\includegraphics[width=0.9\textwidth, angle=0]{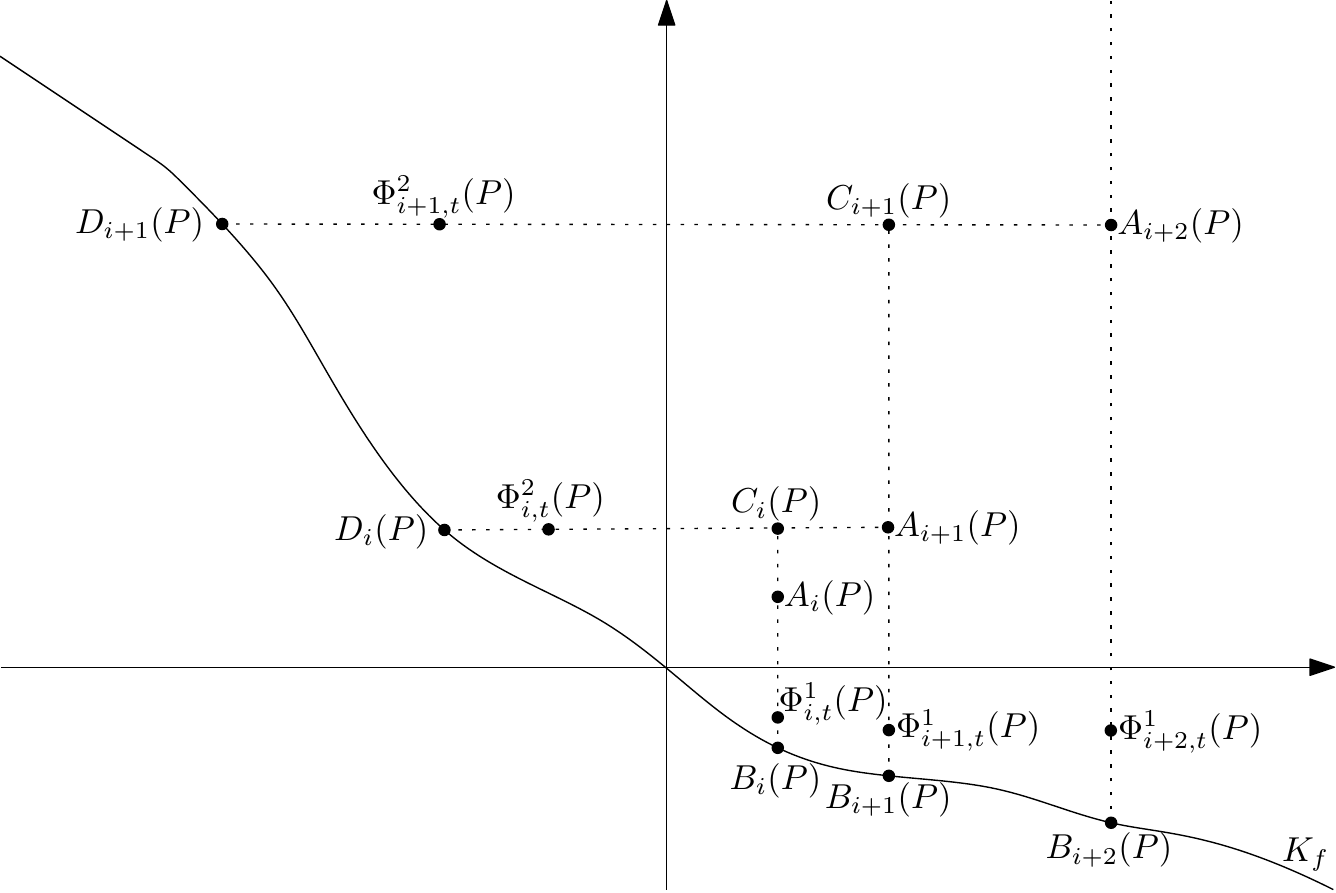}
\centering
\caption{Illustration of the geometric construction of laminates of finite order. This figure illustrates an example of set $K_f$ and the points in $\R^{2 \times 2}_{\sym}$ relevant for our construction. In the figure above, the abscissa represents the $(1,1)$ coordinate in $\R^{2 \times 2}_{\sym}$ and ordinate represents the $(2,2)$ coordinate. The role of the third dimension of $\R^{2 \times 2}_{\sym}$ is ignored. Dotted lines represent rank-one connections.}\label{fig:IllustrationsLaminates}
\end{figure}
For an illustration of the role of the matrices $\Phi_{i,t}^1(P)$ and $\Phi_{i,t}^2(P)$ in the construction, we refer the reader to Figure~\ref{fig:IllustrationsLaminates}. Having defined these maps, the main goal of the remainder of this subsection is to prove the following two lemmas.
\begin{lemma}\label{lemma:MapAiIsOpen}
For all $i \geq 1$, $A_i$ is an open map.
\end{lemma}

\begin{lemma}\label{lemma:InterpolationMapsAreOpen}
For all $t \in (0,1)$ and all $i \geq 1$, $\Phi_{i,t}^1$ and $\Phi_{i,t}^2$ are open maps.
\end{lemma}

\begin{proof}[Proof of Lemma~\ref{lemma:MapAiIsOpen}]
Since the dimension of the domain and the codomain agree, it suffices to show that $A_i$ is a continuous injection (see for instance \cite[Theorem 2B.3]{HAT}). The continuity of the map $A_i$ follows from the fact that $a_i^+$, $a_i^-$, $b$, $\varphi^{\prime}$ are continuous. Now, we only need to show that $A_i$ is injective. Assume that for two $P, \tilde{P} \in \mathcal{P}$, we have
\begin{equation}\label{eq:A_i-equal-for-two-P}
 A_i(P) = A_i(\tilde{P}).
\end{equation}
Write $P = (a_0^+, a_0^-, b)$ and $\tilde{P} = (\tilde{a}_0^+, \tilde{a}_0^-, \tilde{b})$.
From \eqref{eq:A_i-equal-for-two-P}, $b = \tilde{b}$ and $a_i^+(P)= a_i^+(\tilde{P})$ which leads to $a_0^+ = \tilde{a}_0^+$. In addition, \eqref{eq:A_i-equal-for-two-P} also implies that $\varphi^{\prime}(- a_{i+1}^{-}(P)) = \varphi^{\prime}(- a_{i+1}^{-}(\tilde{P}))$ and since $\varphi^{\prime \prime} > 0$, $a_{i+1}^{-}(P) = a_{i+1}^{-}(\tilde{P})$ which leads to $a_0^- = \tilde{a}_0^-$. This proves the injectivity of $A_i$ as wished.
\end{proof}

\begin{proof}[Proof of Lemma~\ref{lemma:InterpolationMapsAreOpen}]
The proof is divided into 2 similar steps. In the first, we treat $\Phi_{i,t}^1$ and the second treats $\Phi_{i,t}^2$. \\
\textbf{Step 1: }We first prove that $\Phi_{i,t}^1$ is an open map for all $t \in (0,1)$. Define the map $\pi_1 \colon \R^{2 \times 2}_{\sym} \to \R^{2 \times 2}_{\sym}$ as 
\[
 \pi_1
 \left(
 \begin{pmatrix}
 a_{11} & a_{12} \\
 a_{12} & a_{22} \\
 \end{pmatrix}
 \right)
 =
 \begin{pmatrix}
 a_{11} & a_{12} \\
 a_{12} & - \varphi^{\prime}(a_{11})
 \end{pmatrix}
\]
and note that $B_i = \pi_1 \circ A_i$ for all $i \geq 1$. Define $\psi_t^1 \colon \R^{2 \times 2}_{\sym} \to \R^{2 \times 2}_{\sym}$ as $\psi_t^1(X) = t \pi_1(X) + (1 - t)X$ and note that $\Phi_{i,t}^1 = \psi_t^1 \circ A_i$. Since $A_i$ is an open map by Lemma~\ref{lemma:MapAiIsOpen}, it suffices to show that $\psi_t^1$ is an open map. To show that $\psi_t^1$ is an open map, we will show that the differential has full rank (see for instance \cite[Proposition 4.28]{LEE}). Instead of considering $\psi_t^1$ as a mapping from $\R^{2 \times 2}_{\sym}$ to $\R^{2 \times 2}_{\sym}$, we choose (without changing the notation) to view $\psi_t^1$ as a mapping from $\R^3$ to $\R^3$. Then 
\[
 \psi_t^1(a_{11}, a_{12}, a_{22}) = (a_{11}, a_{12}, -t \varphi^{\prime}(a_{11}) + (1 - t) a_{22} )
\]
and the differential is
\[
 D \psi_t^1(a_{11}, a_{12}, a_{22}) = 
 \begin{pmatrix}
 1 & 0 & 0 \\
 0 & 1 & 0 \\
 - t \varphi^{\prime \prime}(a_{11}) & 0 & 1-t \\
 \end{pmatrix}.
\]
Then $\det D \psi_t^1(a_{11}, a_{12}, a_{22}) = 1 - t > 0$ for all $t \in (0,1)$ and all $a_{11}, a_{12}, a_{22} \in \R$, which means that $\psi_t^1$ is an open map. Hence $\Phi_{i,t}^1$ is an open map. \\
\textbf{Step 2: }We now prove that $\Phi_{i,t}^2$ is an open map for all $t \in (0,1)$.
Since $\varphi$ is uniformly convex, $\varphi^{\prime}$ admits an inverse map $(\varphi^{\prime})^{-1}$. Define $\pi_2 \colon \R^{2 \times 2}_{\sym} \to \R^{2 \times 2}_{\sym}$ as
\[
 \pi_2
 \left(
 \begin{pmatrix}
 a_{11} & a_{12} \\
 a_{12} & a_{22} \\
 \end{pmatrix}
 \right)
 =
 \begin{pmatrix}
 (\varphi^{\prime})^{-1}(- a_{22}) & a_{12} \\
 a_{12} & a_{22} \\
 \end{pmatrix}.
\]
Note that $D_i = \pi_2 \circ E_i = \pi_2 \circ A_{i+1}$ for all $i \geq 1$. Define $\psi_t^2 \colon \R^{2 \times 2}_{\sym} \to \R^{2 \times 2}_{\sym}$ as $\psi_t^2(X) = t \pi_2(X) + (1-t)X$ and note that $\Phi_{i,t}^2 = \psi_t^2 \circ A_{i+1}$. Since $A_{i+1}$ is an open map, it suffices to show that $\psi_t^2$ is an open map to deduce that $\Phi_{i,t}^2$ is an open map. We can proceed in the same way as in Step 1 to show that $\psi_t^2$ is an open map. This ends the proof.
\end{proof}

\subsection{Construction and parametrisation of interpolation laminates}\label{subsec:ConstructionParametrisationInterpolationLaminates}
In this subsection, the aim is to define new laminates with the help of the interpolation maps from Subsection~\ref{subsec:InterpolationMaps}. As explained in the beginning of this section, these laminates are needed to overcome the difficulty caused by the errors in Lemma~\ref{lemma:ApproxLaminatesFiniteOrder}. Recall the definition of $\mu^{(i)}$ in \eqref{eq:LaminateOfFiniteOrderMu}. We consider this laminate but replace $B_i(P)$ by $\Phi_{i,t}^1(P)$ and $D_i(P)$ by $\Phi_{i,t}^2(P)$ (see Figure~\ref{fig:IllustrationsLaminates}). The coefficients corresponding to $\lambda_i^1(P)$, $\lambda_i^2(P)$ and $\lambda_i^3(P)$ in \eqref{eq:LaminateOfFiniteOrderMu} also have to be adapted. For all $i$ and $t \in (0,1)$, we define
\begin{equation}\label{eq:OneStepInterpolationLaminate}
 \mu^{(i, i+1)}_t(P) \coloneqq \lambda_{i,t}^1(P) \delta_{\Phi_{i,t}^1(P)} + \lambda_{i,t}^2(P) \delta_{\Phi_{i,t}^2(P)} + \lambda_{i,t}^3(P) \delta_{A_{i+1}(P)}
\end{equation}
where
\[
 \lambda_{i,t}^1(P) \coloneqq \frac{\lambda_{i}^1(P)}{t + (1-t)\lambda_{i}^1(P)}, \quad \lambda_{i,t}^2(P) \coloneqq \frac{\lambda_{i}^2(P)}{t + (1-t)\lambda_{i}^1(P)}, \quad \text{and} \quad \lambda_{i,t}^3(P) \coloneqq \frac{t \lambda_{i}^3(P) - (1-t)\lambda_i^2(P) }{t + (1-t)\lambda_{i}^1(P)}.
\]
By Lemma~\ref{lemma:BehaviourOfLambda3}, there is $T_0 \in (0,1)$ and an integer $I_1 \geq 3$ such that for all $i \geq I_1$ and all $t \in (T_0, 1)$ we have
\begin{equation}\label{eq:BehaviourOfLambda3ForInterpolationLaminates}
\frac{1}{r^{2}} < \inf_{t \in (T_0, 1)} \inf_{P \in \mathcal{P}} \lambda_{i,t}^3(P) \leq \sup_{t \in (T_0, 1)} \sup_{P \in \mathcal{P}} \lambda_{i,t}^3(P) < \frac{1}{r^{p + 2 \delta}}.
\end{equation}
By direct computations, the barycenter of $\mu^{(i, i+1)}_t(P)$ is $\overline{\mu^{(i, i+1)}_t(P)} = A_i(P)$. We define a new laminate of finite order $\mu_t^{(i,i+2)}(P)$ by replacing the delta measure $\delta_{A_{i+1}(P)}$ in \eqref{eq:OneStepInterpolationLaminate} by $\mu_t^{(i+1,i+2)}(P)$. In other words,
\begin{align*}
 \mu^{(i, i+2)}_{t}(P) &\coloneqq  \lambda_{i,t}^{1}(P) \delta_{\Phi_{i,t}^1(P)} + \lambda_{i,t}^{2}(P) \delta_{\Phi_{i,t}^2(P)} + \lambda_{i,t}^{3}(P) \lambda_{i+1,t}^{1}(P) \delta_{\Phi_{i+1,t}^1(P)} + \lambda_{i,t}^{3}(P) \lambda_{i+1,t}^{2}(P) \delta_{\Phi_{i+1,t}^2(P)} \\
 & \qquad + \lambda_{i,t}^{3}(P) \lambda_{i+1,t}^{3}(P) \delta_{A_{i+2}(P)}.
\end{align*}
We continue in this way which results in the following definition
\begin{equation}\label{eq:InterpolationLaminate}
\begin{split}
 \mu^{(i,j)}_{t}(P) &\coloneqq \sum_{k = i}^{j-1} \left( \prod_{m = i}^{k-1} \lambda_{m,t}^{3}(P) \right) \lambda_{k,t}^{1}(P) \delta_{\Phi_{k,t}^1(P)} + \sum_{k = i}^{j-1} \left( \prod_{m = i}^{k-1} \lambda_{m,t}^{3}(P) \right) \lambda_{k,t}^{2}(P) \delta_{\Phi_{k,t}^2(P)} \\
 & \qquad + \left( \prod_{m = i}^{j-1} \lambda_{m,t}^{3}(P) \right) \delta_{A_j(P)} \quad \forall i < j, t \in (0,1).
 \end{split}
 \end{equation}
 Note that the barycenter of $\mu_{t}^{(i,j)}(P)$ is $\overline{\mu_{t}^{(i,j)}(P)} = A_i(P)$. 
 In order to correct errors (see Lemma~\ref{lemma:ApproxLaminatesFiniteOrder}) and make sure that the gradients of our sequence of functions progressively approach the set $K_f$, we define the following laminates of finite order
 \begin{equation}\label{eq:CorrectionLaminate1}
 \mu^{1, (i,j)}_{t,t^{\prime}}(P) \coloneqq \dfrac{t}{t^{\prime}} \delta_{\Phi_{i,t^{\prime}}^1(P)} + \left( 1 - \dfrac{t}{t^{\prime}} \right) \mu^{(i,j)}_{t^{\prime}}(P),
\end{equation}
\begin{equation}\label{eq:CorrectionLaminate2}
 \mu^{2, (i,j)}_{t,t^{\prime}}(P) \coloneqq \dfrac{t}{t^{\prime}} \delta_{\Phi_{i,t^{\prime}}^2(P)} + \left( 1 - \dfrac{t}{t^{\prime}} \right) \mu^{(i+1,j)}_{t^{\prime}}(P)
\end{equation}
for all $i < j$, $t, t^{\prime} \in (0,1)$, $t < t^{\prime}$.
Notice that the barycenter of these two laminates are $\overline{\mu^{1, (i,j)}_{t,t^{\prime}}(P)} = \Phi_{i,t}^1(P)$ and $\overline{\mu^{2, (i,j)}_{t,t^{\prime}}(P)} = \Phi_{i,t}^2(P)$. Now, we define collections of sets that will be useful in the proof in the next section. For each $i \geq 1$ and $t \in (0,1)$, we define the sets
\[
 \mathcal{V}_i \coloneqq A_i(\mathcal{P}), \quad \mathcal{U}_i^1 \coloneqq B_i(\mathcal{P}), \quad \mathcal{U}_i^2 \coloneqq D_i(\mathcal{P}), \quad \mathcal{W}_{i,t}^1 \coloneqq \Phi_{i,t}^1(\mathcal{P}), \quad \mathcal{W}_{i,t}^2 \coloneqq \Phi_{i,t}^2(\mathcal{P}).
\]
By Lemma~\ref{lemma:MapAiIsOpen} and \ref{lemma:InterpolationMapsAreOpen}, the sets $\mathcal{V}_i$, $\mathcal{W}_{i,t}^1$ and $\mathcal{W}_{i,t}^2$ are open for all $i \geq 1$ and $t \in (0,1)$. The laminates \eqref{eq:OneStepInterpolationLaminate}, \eqref{eq:CorrectionLaminate1}, \eqref{eq:CorrectionLaminate2}, the openness of the sets $\mathcal{V}_i$, $\mathcal{W}_{i,t}^1$ and $\mathcal{W}_{i,t}^2$ combined with Lemma~\ref{lemma:ApproxLaminatesFiniteOrder} (see also Remark~\ref{rmk:AboutHowBuildingBlocksAreUsed}) immediately leads to the three propositions that follow below. We shall provide an explicit proof only for Proposition~\ref{prop:CVCorr1}, because the proof of Proposition~\ref{prop:CVCorr2} is essentially the same and the proof of Proposition~\ref{prop:CVIteration} is in the same spirit but slightly less involved. Recall the constants $I_1$ and $T_0$ defined in \eqref{eq:BehaviourOfLambda3ForInterpolationLaminates}.

\begin{proposition}\label{prop:CVIteration}
 Let $j > i \geq I_1$, $P \in \mathcal{P}$, $\alpha \in (0,1)$, $t \in (T_0, 1)$, $b \in \R^2$ and $U$ convex and open. For any $\eps > 0$, there is a piecewise affine Lipschitz mapping $w \colon U \to \R^2$ such that
 \begin{enumerate}
 \item $\| w - (A_i(P)x + b) \|_{C^{\alpha}(\overline{U}, \R^2)} < \eps$;
 \smallskip
 \item $w = A_i(P)x + b$ on $\partial U$;
 \smallskip
 \item $Dw(x) \in \bigcup_{k = i}^{j-1} \mathcal{W}_{k,t}^1 \cup \bigcup_{k = i}^{j-1} \mathcal{W}_{k,t}^2 \cup \mathcal{V}_{j}$ for a.e. $x \in U$;
 \smallskip
 \item $(1 - \frac{1}{r^p}) ( \frac{1}{r^2} )^{k-i} |U| \leq |\{ Dw \in \mathcal{W}_{k, t}^1 \cup \mathcal{W}_{k, t}^2 \}| \leq ( 1 - \frac{1}{r^2}) ( \frac{1}{r^{p + 2 \delta}} )^{k-i} |U|$ for all $i \leq k < j$;
\smallskip
 \item $(\frac{1}{r^2})^{j-i} |U| \leq |\{ Dw \in \mathcal{V}_{j} \}| \leq (\frac{1}{r^{p + 2 \delta}})^{j-i} |U| $.
 \end{enumerate}
\end{proposition}
\begin{proposition}\label{prop:CVCorr1}
 Let $j > i \geq I_1$, $P \in \mathcal{P}$, $\alpha \in (0,1)$, $T_0 < t < t^{\prime} < 1$, $b \in \R^2$ and $U$ convex and open. For any $\eps > 0$, there is a piecewise affine Lipschitz mapping $w \colon U \to \R^2$ such that
 \begin{enumerate}
 \item $\| w - (\Phi_{i,t}^1(P) x + b) \|_{C^{\alpha}(\overline{U}, \R^2)} < \eps$;
 \smallskip
 \item $w = \Phi_{i,t}^1(P) x + b$ on $\partial U$;
 \smallskip
 \item $Dw(x) \in \bigcup_{k = i}^{j-1} \mathcal{W}_{k,t^{\prime}}^1 \cup \bigcup_{k = i}^{j-1} \mathcal{W}_{k,t^{\prime}}^2 \cup \mathcal{V}_{j}$ for a.e. $x \in U$;
 \smallskip
 \item $\frac{t}{t^{\prime}} |U| \leq |\{ Dw \in \mathcal{W}_{i, t^{\prime}}^1 \cup \mathcal{W}_{i, t^{\prime}}^2 \}| \leq |U|$;
\smallskip
 \item $(1 - \frac{1}{r^p}) (1 - \frac{t}{t^{\prime}}) ( \frac{1}{r^2} )^{k-i} |U| \leq |\{ Dw \in \mathcal{W}_{k, t^{\prime}}^1 \cup \mathcal{W}_{k, t^{\prime}}^2 \}| \, \leq \, ( 1 - \frac{1}{r^2}) (1 - \frac{t}{t^{\prime}}) ( \frac{1}{r^{p + 2 \delta}} )^{k-i} |U|$ for all $i < k < j$;
\smallskip
 \item $(1 - \frac{t}{t^{\prime}}) (\frac{1}{r^2})^{j-i} |U| \leq |\{ Dw \in \mathcal{V}_{j} \}| \leq (1 - \frac{t}{t^{\prime}}) (\frac{1}{r^{p + 2 \delta}})^{j-i} |U| $.
 \end{enumerate}
\end{proposition}
\begin{proposition}\label{prop:CVCorr2}
 Let $j > i \geq I_1$, $P \in \mathcal{P}$, $\alpha \in (0,1)$, $T_0 < t < t^{\prime} < 1$, $b \in \R^2$ and $U$ convex and open. For any $\eps > 0$, there is a piecewise affine Lipschitz mapping $w \colon U \to \R^2$ such that
 \begin{enumerate}
 \item $\| w - (\Phi_{i,t}^2(P) x + b) \|_{C^{\alpha}(\overline{U}, \R^2)} < \eps$;
 \smallskip
 \item $w = \Phi_{i,t}^2(P) x + b$ on $\partial U$;
 \smallskip
 \item $Dw(x) \in \bigcup_{k = i}^{j-1} \mathcal{W}_{k,t^{\prime}}^1 \cup \bigcup_{k = i}^{j-1} \mathcal{W}_{k,t^{\prime}}^2 \cup \mathcal{V}_{j}$ for a.e. $x \in U$;
 \smallskip
 \item $\frac{t}{t^{\prime}} |U| \leq |\{ Dw \in \mathcal{W}_{i, t^{\prime}}^1 \cup \mathcal{W}_{i, t^{\prime}}^2 \}| \leq |U|$;
\smallskip
 \item $(1 - \frac{1}{r^p}) (1 - \frac{t}{t^{\prime}}) ( \frac{1}{r^2} )^{k-i-1} |U| \leq |\{ Dw \in \mathcal{W}_{k, t^{\prime}}^1 \cup \mathcal{W}_{k, t^{\prime}}^2 \}| \leq ( 1 - \frac{1}{r^2}) (1 - \frac{t}{t^{\prime}}) ( \frac{1}{r^{p + 2 \delta}} )^{k-i-1} |U|$ for all $i < k < j$;
\smallskip
 \item $(1 - \frac{t}{t^{\prime}}) (\frac{1}{r^2})^{j-i-1} |U| \leq |\{ Dw \in \mathcal{V}_{j} \}| \leq (1 - \frac{t}{t^{\prime}}) (\frac{1}{r^{p + 2 \delta}})^{j-i-1} |U| $.
 \end{enumerate}
\end{proposition}

\begin{proof}[Proof of Proposition~\ref{prop:CVCorr1}]
The laminate of finite order $ \mu^{1, (i,j)}_{t,t^{\prime}}(P)$ from \eqref{eq:CorrectionLaminate1} satisfies the property that $\overline{ \mu^{1, (i,j)}_{t,t^{\prime}}(P)} = \Phi_{i,t}^1(P)$. Moreover, the following estimates follow from the definition of $ \mu^{1, (i,j)}_{t,t^{\prime}}(P)$:
\begin{align*}
 & \mu^{1, (i,j)}_{t,t^{\prime}}(P)\left( \bigcup_{k = i}^{j-1} \mathcal{W}_{k,t^{\prime}}^1 \cup \bigcup_{k = i}^{j-1} \mathcal{W}_{k,t^{\prime}}^2 \cup \mathcal{V}_{j} \right) = 1;\\
 &\frac{t}{t^{\prime}} \leq  \mu^{1, (i,j)}_{t,t^{\prime}}(P)\left( \mathcal{W}_{i,t^{\prime}}^1 \cup \mathcal{W}_{i,t^{\prime}}^2 \right) \leq 1;\\
 &\left( 1 - \frac{1}{r^p} \right) \left( 1 - \frac{t}{t^{\prime}} \right) \left( \frac{1}{r^2} \right)^{k-i}  \leq  \mu^{1, (i,j)}_{t,t^{\prime}}(P)\left( \mathcal{W}_{k,t^{\prime}}^1 \cup \mathcal{W}_{k,t^{\prime}}^2 \right) \leq \left( 1 - \frac{1}{r^2} \right) \left( 1 - \frac{t}{t^{\prime}} \right) \left( \frac{1}{r^{p + 2 \delta}} \right)^{k-i} \, \forall i < k < j;\\
 &\left( 1 - \frac{t}{t^{\prime}} \right) \left( \frac{1}{r^{2}} \right)^{j-i}  \leq  \mu^{1, (i,j)}_{t,t^{\prime}}(P) \left( \mathcal{V}_{j} \right) \leq \left( 1 - \frac{t}{t^{\prime}} \right) \left( \frac{1}{r^{p + 2 \delta}} \right)^{j-i}.\\
\end{align*}
Therefore, since the sets $\mathcal{W}_{k,t^{\prime}}^1$, $\mathcal{W}_{k,t^{\prime}}^2$ ($i < k < j$) and $\mathcal{V}_j$ are open, by applying Lemma~\ref{lemma:ApproxLaminatesFiniteOrder} (see also Remark~\ref{rmk:AboutHowBuildingBlocksAreUsed}), we find a piecewise affine Lipschitz mapping $w \colon U \to \R^2$ satisfying all the desired properties. This concludes the proof of Proposition~\ref{prop:CVCorr1}.
\end{proof}

Finally, we end this subsection and section by the following observation: 
There is a finite constant $C > 1$ such that
\begin{equation}\label{eq:ExpontialBoundInTheSets}
|X_{11}|^p + |X_{22}|^p < C r^{p \ell}, \, |X_{12}| < 1, \, |X_{21}| < 1 \quad \forall X \in \overline{\mathcal{U}}_{\ell}^1 \cup \overline{\mathcal{U}}_{\ell}^2 \cup \mathcal{V}_{\ell} \cup \bigcup_{t \in (T_0,1), r = 1,2} \mathcal{W}_{\ell,t}^r.
\end{equation}
This follows from the definition of the sets $\mathcal{U}_{\ell}^1$, $\mathcal{U}_{\ell}^2$, $\mathcal{V}_{\ell}$, $\mathcal{W}_{\ell, t}^1$ and $\mathcal{W}_{\ell, t}^2$ combined with direct computations.

\section{Convex integration scheme (Part 2 of the proof of Theorem~\ref{thm:MainFunctionals})}\label{sec:Scheme}
In this section, we describe the convex integration scheme and finish the proof of Theorem~\ref{thm:MainFunctionals}. More precisely, we describe the scheme in Subsection~\ref{subsec:DescriptionScheme}. At the end of this first subsection, we will have finished the proof of Theorem~\ref{thm:MainFunctionals}. In order to do this, we refer to several properties that we prove in the remaining subsections.
\subsection{Description of the Convex Integration Scheme}\label{subsec:DescriptionScheme}
Recall $I_1$ and $T_0$ defined in \eqref{eq:BehaviourOfLambda3ForInterpolationLaminates}.
Select\footnote{To select such $I$ and $T_0^{\prime}$, it suffices to choose $T_0^{\prime}$ and $\tilde{I} \geq I_1$ such that for all $i \geq \tilde{I}$, $\mathcal{W}_{i, T_0^{\prime}}^1, \mathcal{W}_{i, T_0^{\prime}}^2 \subset \{ X \in \R^{2 \times 2}_{\sym} : X_{22} < 0 \}$. This is possible because $0 < \lambda, \Lambda$. Hence \eqref{eq:AboutDistancesBetweenSets} holds for any $I \geq \tilde{I}$ and $T_0^{\prime}$. Then it only remains to select $I \geq \tilde{I}$ so that \eqref{eq:AboutIToThePowerTwo} holds.} $I \geq I_1$ and $0 < T_0^{\prime} < 1$ such that 
\begin{align}
 r^I &\geq \frac{1}{r^{2 \delta} - 1}; \label{eq:AboutIToThePowerTwo}\\
1 &< \operatorname{dist}\left({\bigcup_{j = 0}^{\infty} \overline{\mathcal{V}}_{I + j}},{\bigcup_{\substack{j \geq 0, \\ t > T^{\prime}_0}} \overline{\mathcal{W}}_{I + j, t}^1 \cup \bigcup_{\substack{j \geq 0, \\ t > T^{\prime}_0}} \overline{\mathcal{W}}_{I + j, t}^2 }\right).
\label{eq:AboutDistancesBetweenSets}
\end{align}
Take $M \in \mathcal{V}_I$. Define $v_I \colon \Omega \to \R^2$ as $v_I(x) \coloneqq Mx$. Let $\{ \delta_{\ell} \}_{\ell = 1}^{\infty}$ be a decreasing sequence such that $\delta_{\ell} \leq 2 \delta$ for all $\ell$ and $\delta_{\ell} \to \delta$. Let $\{ t_{\ell} \}_{\ell = 1}^{\infty}$ be an increasing sequence such that (recall $C_{r,p}$ defined in \eqref{eq:DefinitionOfCRP}):
\begin{align}
 1 - \frac{1}{r^{I+2 \ell}} &\leq t_{\ell} \leq 1; \label{eq:SeqLowerBound}\\
 2 C_{r,p} (t_{\ell + 1} - t_{\ell}) &\leq \left( \frac{1}{r^{p + \delta_{\ell + 1}}} \right)^{j} - \left( \frac{1}{r^{p + \delta_{\ell}}} \right)^{j} \, \forall j = 1, \ldots, \ell + 1; \label{eq:Difficult} \\
 t_{1} &\geq \max \left(\frac{9}{10}, T_0, T_0^{\prime} \right). \label{eq:Easy}
\end{align}
The existence of the sequences $\{ \delta_{\ell} \}_{\ell = 1}^{\infty}$ and $\{ t_{\ell} \}_{\ell = 1}^{\infty}$ will be proved in Subsection~\ref{subsec:tlseq}. Let $\eps > 0$ be arbitrary. 

\noindent \textbf{Step 1:}
We describe the first step of the scheme. 
Let $0 < \eps_I < \sfrac{\eps}{2}$.
Since $M \in \mathcal{V}_I$ there is $P \in \mathcal{P}$ such that $M = A_I(P)$. By Proposition~\ref{prop:CVIteration}, there is a piecewise affine mapping $v_{I+1} \colon \Omega \to \R^2$ such that
\begin{enumerate}
\item $\| v_{I+1} - v_I \|_{C^{\alpha}(\overline{\Omega}; \R^2)} < \frac{\eps_{I}}{2^{I}}$; \label{item:CIS_A1}
\smallskip
\item $Dv_{I+1}(x) \in \mathcal{W}_{I,t}^1 \cup \mathcal{W}_{I,t}^2 \cup \mathcal{V}_{I+1}$ for a.e. $x \in \Omega$; \label{item:CIS_A2}
\smallskip
\item $(1 - \frac{1}{r^p}) |\Omega| < |\{ Dv_{I+1} \in \mathcal{W}_{I, t_1}^1 \cup \mathcal{W}_{I, t_1}^2 \}| < ( 1 - \frac{1}{r^2}) |\Omega|$; \label{item:CIS_A3}
\smallskip
\item $\frac{1}{r^2} |\Omega| \leq |\{ Dv_{I+1} \in \mathcal{V}_{I+1} \}| \leq \frac{1}{r^{p + 2 \delta}} |\Omega| \leq \frac{1}{r^{p + \delta_{1}}} |\Omega|$. \label{item:CIS_A4}
\end{enumerate}
Let $\mathcal{F}_{I+1}$ be a collection of open convex mutually disjoint sets covering $\Omega$ up to a set of measure 0 so that $v_{I+1} |_{U}$ is affine for each $U \in \mathcal{F}_{I+1}$. Without loss of generality, $\sup_{U \in \mathcal{F}_{I+1}} \diam{U} < 2^{-(I+1)}$. Select $0 < \eps_{I+1} < \min(\eps_I, 2^{-(I+1)})$ such that $\| Dv_{I+1} \ast \rho_{\eps_{I+1}} - Dv_{I+1} \|_{L^1({\Omega, \R^{2 \times 2}})} < 2^{-(I+1)}$.

\noindent\textbf{Inductive step, $\pmb{(\ell+1)}$-th step: } We assume that we are given $v_{I + \ell} \colon \Omega \to \R^2$, $v_{I + \ell - 1} \colon \Omega \to \R^2$, $\eps_{I + \ell - 1}$ and $\eps_{I + \ell}$ such that
\begin{enumerate}
\setcounter{enumi}{4}
\item $\| v_{I+\ell} - v_{I+\ell-1} \|_{C^{\alpha}(\overline{\Omega}; \R^2)} < \frac{\eps_{I+\ell-1}}{2^{I + \ell-1}}$; \label{item:CIS_G1}
\smallskip
\item $Dv_{I + \ell}(x) \in \bigcup_{j = 0}^{\ell - 1} \mathcal{W}_{I+j,t_\ell}^1 \cup \bigcup_{j = 0}^{\ell - 1} \mathcal{W}_{I+j,t_\ell}^2 \cup \mathcal{V}_{I+\ell}$ for a.e. $x \in \Omega$; \label{item:CIS_G2}
\smallskip
\item $(1 - \frac{1}{r^p}) ( \frac{1}{r^2} )^j \frac{t_{j+1}}{t_{\ell}} |\Omega| \leq |\{ Dv_{I + \ell} \in \mathcal{W}_{I+j, t_\ell}^1 \cup \mathcal{W}_{I+j, t_\ell}^2 \}| \leq ( 1 - \frac{1}{r^2}) ( \frac{1}{r^{p + \delta_\ell}} )^j |\Omega|$ for all $0 \leq j < \ell$; \label{item:CIS_G3}
\smallskip
 \item $(\frac{1}{r^2})^\ell |\Omega| \leq |\{ Dv_{I + \ell} \in \mathcal{V}_{I+\ell} \}| \leq (\frac{1}{r^{p + \delta_\ell}})^\ell |\Omega|$; \label{item:CIS_G4}
 \smallskip
\item $\| Dv_{I+\ell} \ast \rho_{\eps_{I+\ell}} - Dv_{I+\ell} \|_{L^1({\Omega, \R^{2 \times 2}})} < 2^{-(I+\ell)}$; \label{item:CIS_G5}
 \smallskip
\item there is a collection $\mathcal{F}_{I+\ell}$ of open convex mutually disjoint sets covering $\Omega$ up to set of measure 0 such that $v_{I+\ell} |_U$ is affine for all $U \in \mathcal{F}_{I+\ell}$ and $\sup_{U \in \mathcal{F}_{I+\ell}} \diam{U} < 2^{-(I+\ell)}$. \label{item:CIS_G6}
\end{enumerate}
We will build $v_{I + \ell + 1}$ by modifying $v_{I+\ell}$ on each $U \in \mathcal{F_{I + \ell}}$.
For each $U \in \mathcal{F}_{I+\ell}$, the construction is as follows:
\\
\fbox{If $Dv_{I+\ell} |_U \in \mathcal{W}_{I+j,t_\ell}^1$ for some $0 \leq j < \ell$: }
There is $P \in \mathcal{P}$ such that $Dv_{I+\ell} |_U = \Phi_{I + j, t_\ell}^1(P)$.
By Proposition~\ref{prop:CVCorr1}, there is a piecewise affine $v_{I + \ell +1, U} \colon U \to \R^2$ such that
\begin{enumerate}
\setcounter{enumi}{10}
\item $\| v_{I+\ell+1, U} - v_{I+\ell} \|_{C^{\alpha}(\overline{U}; \R^2)} < \frac{\eps_{I+\ell}}{2^{I + \ell + 1}}$; \label{item:CIS_H1}
\smallskip
 \item $v_{I + \ell + 1, U} = v_{I + \ell}$ on $\partial U$;
\smallskip
\item $Dv_{I+\ell+1, U}(x) \in \bigcup_{k = j}^{\ell} \mathcal{W}_{I+k,t_{\ell+1}}^1 \cup \bigcup_{k = j}^{\ell} \mathcal{W}_{I+k,t_{\ell+1}}^2 \cup \mathcal{V}_{I+\ell+1}$ for a.e. $x \in \Omega$; \label{item:CIS_H2}
\smallskip
\item $\frac{t_{\ell}}{t_{\ell + 1}} |U| \leq |\{ Dv_{I+\ell+1, U} \in \mathcal{W}_{I+j, t_{\ell + 1}}^1 \cup \mathcal{W}_{I+j, t_{\ell + 1}}^2 \}| \leq |U|$; \label{item:CIS_H3}
\smallskip
\item $( 1 - \frac{t_{\ell}}{t_{\ell + 1}}) (\frac{1}{r^2})^m (1 - \frac{1}{r^p}) |U| \leq |\{ Dv_{I+\ell, U} \in \mathcal{W}_{I+j+m, t_{\ell + 1}}^1 \cup \mathcal{W}_{I+j+m, t_{\ell + 1}}^2 \}| \newline \leq ( 1 - \frac{t_{\ell}}{t_{\ell + 1}}) {(\frac{1}{r^{p + 2 \delta}})^m ( 1 - \frac{1}{r^2}) |U|}$ for all $1 \leq m \leq \ell - j$; \label{item:CIS_H4}
\smallskip
 \item $( 1 - \frac{t_{\ell}}{t_{\ell + 1}}) (\frac{1}{r^2})^{\ell + 1 - j} |U| \leq |\{ Dv_{I+\ell+1, U} \in \mathcal{V}_{I+\ell+1} \}| \leq ( 1 - \frac{t_{\ell}}{t_{\ell + 1}}) (\frac{1}{r^{p + 2 \delta}})^{\ell + 1 - j} |U| $. \label{item:CIS_H5}
\end{enumerate}
\fbox{If $Dv_{I+\ell} |_U \in \mathcal{W}_{I+j,t_\ell}^2$ for some $0 \leq j < \ell$: }
There is $P \in \mathcal{P}$ such that $Dv_{I+\ell} |_U = \Phi_{I + j, t_\ell}^2(P)$.
By Proposition~\ref{prop:CVCorr2}, there is a piecewise affine $v_{I + \ell +1, U} \colon U \to \R^2$ such that
\begin{enumerate}
\setcounter{enumi}{16}
\item $\| v_{I+\ell+1, U} - v_{I+\ell} \|_{C^{\alpha}(\overline{U}; \R^2)} < \frac{\eps_{I+\ell}}{2^{I + \ell + 1}}$; \label{item:CIS_I1}
\smallskip
 \item $v_{I + \ell + 1, U} = v_{I + \ell}$ on $\partial U$;
\smallskip
\item $Dv_{I+\ell+1, U}(x) \in \bigcup_{k = j+1}^{\ell} \mathcal{W}_{I+k,t_{\ell+1}}^1 \cup \bigcup_{k = j}^{\ell} \mathcal{W}_{I+k,t_{\ell+1}}^2 \cup \mathcal{V}_{I+\ell+1}$ for a.e. $x \in \Omega$; \label{item:CIS_I2}
\smallskip
\item $\frac{t_{\ell}}{t_{\ell + 1}} |U| \leq |\{ Dv_{I+\ell+1, U} \in \mathcal{W}_{I+j, t_{\ell + 1}}^1 \cup \mathcal{W}_{I+j, t_{\ell + 1}}^2 \}| \leq |U|$; \label{item:CIS_I3}
\smallskip
\item $( 1 - \frac{t_{\ell}}{t_{\ell + 1}}) (\frac{1}{r^2})^{m-1} (1 - \frac{1}{r^p}) |U| \leq |\{ Dv_{I+\ell, U} \in \mathcal{W}_{I+j+m, t_{\ell +1}}^1 \cup \mathcal{W}_{I+j+m, t_{\ell + 1}}^2 \}| \newline \leq ( 1 - \frac{t_{\ell}}{t_{\ell + 1}}) (\frac{1}{r^{p + 2 \delta}})^{m-1} ( 1 - \frac{1}{r^2}) |U|$ for all $1 \leq m \leq \ell - j$; \label{item:CIS_I4}
\smallskip
 \item $( 1 - \frac{t_{\ell}}{t_{\ell + 1}}) (\frac{1}{r^2})^{\ell - j} |U| \leq |\{ Dv_{I+\ell+1, U} \in \mathcal{V}_{I+\ell+1} \}| \leq ( 1 - \frac{t_{\ell}}{t_{\ell + 1}}) (\frac{1}{r^{p + 2 \delta}})^{\ell - j} |U| $. \label{item:CIS_I5}
\end{enumerate}
\fbox{If $Dv_{I+\ell} |_U \in \mathcal{V}_{I+\ell}$: }
There is $P \in \mathcal{P}$ such that $Dv_{I+\ell} |_U = A_{I+\ell}(P)$.
By Proposition~\ref{prop:CVIteration}, there is a piecewise affine $v_{I + \ell +1, U} \colon U \to \R^2$ such that
\begin{enumerate}
\setcounter{enumi}{22}
\item $\| v_{I+\ell+1, U} - v_{I+\ell} \|_{C^{\alpha}(\overline{U}; \R^2)} < \frac{\eps_{I+\ell}}{2^{I + \ell + 1}}$; \label{item:CIS_J1}
\smallskip
 \item $v_{I + \ell + 1, U} = v_{I + \ell}$ on $\partial U$;
\smallskip
\item $Dv_{I+\ell+1, U}(x) \in \mathcal{W}_{I+\ell, t_{\ell+1}}^1 \cup \mathcal{W}_{I+\ell,t_{\ell+1}}^2 \cup \mathcal{V}_{I+\ell+1}$ for a.e. $x \in \Omega$; \label{item:CIS_J2}
\smallskip
\item $(1 - \frac{1}{r^p}) |U| \leq |\{ Dv_{I+\ell+1, U} \in \mathcal{W}_{I+\ell, t_{\ell + 1}}^1 \cup \mathcal{W}_{I+\ell, t_{\ell + 1}}^2 \}| \leq ( 1 - \frac{1}{r^2}) |U|$; \label{item:CIS_J3}
 \item $\frac{1}{r^2} |U| \leq |\{ Dv_{I+\ell+1, U} \in \mathcal{V}_{I+\ell+1} \}| \leq \frac{1}{r^{p + 2 \delta}} |U|$. \label{item:CIS_J4}
\end{enumerate}
Now, we can define $v_{I + \ell + 1}$. We define the function $v_{I + \ell + 1} \colon \Omega \to \R^2$ as $v_{I + \ell + 1}(x) = v_{I + \ell + 1, U}(x)$ if $x \in U$ for all $U \in \mathcal{F}_{I + \ell}$. 
Now we estimate the following quantities (for which we also introduce new notation):
\begin{itemize}
 \item $\mathcal{Q}^{\text{Step $(\ell+1)$}}_{I+j} \coloneqq |\{ Dv_{I+\ell+1} \in \mathcal{W}_{I+j, t_{\ell+1}}^1 \cup \mathcal{W}_{I+j, t_{\ell+1}}^2 \}|$  for all $0 \leq j \leq \ell$;
 \item $\widetilde{\mathcal{Q}}^{\text{Step $(\ell+1)$}} \coloneqq |\{ Dv_{I+\ell+1} \in \mathcal{V}_{I+\ell+1} \}|$.
\end{itemize}
For $j = 0$, we find
\begin{align*}
\mathcal{Q}^{\text{Step $(\ell+1)$}}_{I} \overset{\ref{item:CIS_H3}, \ref{item:CIS_I3}}&{\leq} |\{ Dv_{I+\ell} \in \mathcal{W}_{I, t_{\ell}}^1 \}| + |\{ Dv_{I+\ell} \in \mathcal{W}_{I, t_{\ell}}^2 \}| \\
&\leq |\{ Dv_{I+\ell} \in \mathcal{W}_{I, t_{\ell}}^1 \cup \mathcal{W}_{I, t_{\ell}}^2 \}|\overset{\ref{item:CIS_G3}}{\leq} \left( 1 - \frac{1}{r^2} \right) |\Omega|; \\
\mathcal{Q}^{\text{Step $(\ell+1)$}}_{I} \overset{\ref{item:CIS_H3}, \ref{item:CIS_I3}}&{\geq} \frac{t_{\ell}}{t_{\ell + 1}} |\{ Dv_{I+\ell} \in \mathcal{W}_{I, t_{\ell}}^1 \cup \mathcal{W}_{I, t_{\ell}}^2 \}|  \\
\overset{\ref{item:CIS_G3}}&{\geq} \frac{t_{\ell}}{t_{\ell + 1}} \left( 1 - \frac{1}{r^p} \right) \frac{t_{1}}{t_{\ell}} |\Omega| \geq \left( 1 - \frac{1}{r^p} \right) \frac{t_{1}}{t_{\ell+1}} |\Omega|.
\end{align*}
For $1 \leq j < \ell$, we find
\begin{align*}
\mathcal{Q}^{\text{Step $(\ell+1)$}}_{I+j} \overset{\ref{item:CIS_H3}, \ref{item:CIS_H4}, \ref{item:CIS_I3}, \ref{item:CIS_I4}}&{\leq} \sum_{k = 0}^{j-1}  \left( 1 - \frac{t_{\ell}}{t_{\ell + 1}} \right) \left( \frac{1}{r^{p + 2 \delta}} \right)^{j-k} \left( 1 - \frac{1}{r^2} \right) |\{ Dv_{I+\ell} \in \mathcal{W}_{I+k, t_{\ell}}^1 \}| \\
&\qquad + \sum_{k = 0}^{j-1}  \left( 1 - \frac{t_{\ell}}{t_{\ell + 1}} \right) \left( \frac{1}{r^{p + 2 \delta}} \right)^{j-k-1} \left( 1 - \frac{1}{r^2} \right) |\{ Dv_{I+\ell} \in \mathcal{W}_{I+k, t_{\ell}}^2 \}| \\
&\qquad + |\{ Dv_{I+\ell} \in \mathcal{W}_{I+j, t_{\ell}}^1 \}| + |\{ Dv_{I+\ell} \in \mathcal{W}_{I+j, t_{\ell}}^2 \}| \\
&\leq \left( 1 - \frac{1}{r^2} \right) \left( 1 - \frac{t_{\ell}}{t_{\ell + 1}} \right) \sum_{k = 0}^{j-1} \left( \frac{1}{r^{p + \delta_1}} \right)^{j-k-1} |\{ Dv_{I+\ell} \in \mathcal{W}_{I+k, t_{\ell}}^1 \cup \mathcal{W}_{I+k, t_{\ell}}^2 \}| \\
&\qquad + |\{ Dv_{I+\ell} \in \mathcal{W}_{I+j, t_{\ell}}^1 \cup \mathcal{W}_{I+j, t_{\ell}}^2 \}| \\
\overset{\ref{item:CIS_G3}}&{\leq} \left( 1 - \frac{1}{r^2} \right)^2 \left( 1 - \frac{t_{\ell}}{t_{\ell + 1}} \right) \left( \frac{1}{r^{p + \delta_1}} \right)^{j-1} |\Omega| \sum_{k = 0}^{j-1} \left( \frac{r^{p + \delta_1}}{r^{p + \delta_\ell}} \right)^k + \left( 1 - \frac{1}{r^2} \right) \left( \frac{1}{r^{p + \delta_\ell}} \right)^j |\Omega| \\
&\leq \left( 1 - \frac{1}{r^2} \right)^2 \left(1 - \frac{t_{\ell}}{t_{\ell + 1}} \right) j \left(\frac{r^{(\delta_1 - \delta_{\ell})}}{r^{p + \delta_1}} \right)^{j-1} |\Omega| + \left( 1 - \frac{1}{r^2} \right) \left( \frac{1}{r^{p + \delta_\ell}} \right)^j |\Omega| \\
\overset{\eqref{eq:DefinitionOfCRP}, \eqref{eq:Easy}}&{\leq} \left( 1 - \frac{1}{r^2} \right) \left( 2 C_{r,p}\left( 1 - \frac{1}{r^2} \right) (t_{\ell+1} - t_\ell) + \left( \frac{1}{r^{p + \delta_\ell}} \right)^j \right) |\Omega| \\
&\leq \left( 1 - \frac{1}{r^2} \right) \left( 2 C_{r,p} (t_{\ell+1} - t_\ell) + \left( \frac{1}{r^{p + \delta_\ell}} \right)^j \right) |\Omega|
\overset{\eqref{eq:Difficult}}{\leq} \left( 1 - \frac{1}{r^2} \right) \left( \frac{1}{r^{p + \delta_{\ell + 1}}} \right)^j |\Omega|; \\
\mathcal{Q}^{\text{Step $(\ell+1)$}}_{I+j} \overset{\ref{item:CIS_H3}, \ref{item:CIS_I3}}&{\geq} \frac{t_{\ell}}{t_{\ell + 1}} |\{ Dv_{I+\ell} \in \mathcal{W}_{I+j, t_{\ell}}^1 \cup \mathcal{W}_{I+j, t_{\ell}}^2 \}|  \\
\overset{\ref{item:CIS_G3}}&{\geq} \frac{t_{\ell}}{t_{\ell + 1}} \left( 1 - \frac{1}{r^p} \right) \left( \frac{1}{r^2} \right)^j \frac{t_{j+1}}{t_{\ell}} |\Omega| \geq \left( 1 - \frac{1}{r^p} \right) \left( \frac{1}{r^2} \right)^j \frac{t_{j+1}}{t_{\ell+1}} |\Omega|.
\end{align*}
For $j = \ell + 1$, we find
\begin{align*}
\mathcal{Q}^{\text{Step $(\ell + 1)$}}_{I+\ell} \overset{\ref{item:CIS_H3}, \ref{item:CIS_I3}, \ref{item:CIS_J3}}&{\leq} \sum_{k = 0}^{\ell-1}  \left( 1 - \frac{t_{\ell}}{t_{\ell + 1}} \right) \left( \frac{1}{r^{p + 2 \delta}} \right)^{\ell-k} \left( 1 - \frac{1}{r^2} \right) |\{ Dv_{I+\ell} \in \mathcal{W}_{I+k, t_{\ell}}^1 \}| \\
&\qquad + \sum_{k = 0}^{\ell-1}  \left(1 - \frac{t_{\ell}}{t_{\ell + 1}} \right) \left( \frac{1}{r^{p + 2 \delta}} \right)^{\ell-k-1} \left( 1 - \frac{1}{r^2} \right) |\{ Dv_{I+\ell} \in \mathcal{W}_{I+k, t_{\ell}}^2 \}| \\
& \qquad + \left( 1 - \frac{1}{r^2} \right) |\{ Dv_{I+\ell} \in \mathcal{V}_{I+\ell} \}| \\
&\leq \left( 1 - \frac{1}{r^2} \right) \left( 1 - \frac{t_{\ell}}{t_{\ell + 1}} \right) \sum_{k = 0}^{\ell-1} \left( \frac{1}{r^{p + \delta_1}} \right)^{\ell-k-1} |\{ Dv_{I+\ell} \in \mathcal{W}_{I+k, t_{\ell}}^1 \cup \mathcal{W}_{I+k, t_{\ell}}^2 \}| \\
& \qquad + \left( 1 - \frac{1}{r^2} \right) |\{ Dv_{I+\ell} \in \mathcal{V}_{I+\ell} \}| \\
\overset{\ref{item:CIS_G3}, \ref{item:CIS_G4}}&{\leq} \left( 1 - \frac{1}{r^2} \right)^2 \left( 1 - \frac{t_{\ell}}{t_{\ell + 1}} \right) \left( \frac{1}{r^{p + \delta_1}} \right)^{\ell-1} |\Omega| \sum_{k = 0}^{\ell-1} \left( \frac{r^{p + \delta_1}}{r^{p + \delta_\ell}} \right)^k + \left( 1 - \frac{1}{r^2} \right) \left( \frac{1}{r^{p + \delta_\ell}} \right)^\ell |\Omega| \\
 &\leq \left( 1 - \frac{1}{r^2} \right)^2 \left( 1 - \frac{t_{\ell}}{t_{\ell + 1}} \right) \ell \left( \frac{r^{(\delta_1 - \delta_{\ell})}}{r^{p + \delta_1}} \right)^{\ell - 1} |\Omega| + \left( 1 - \frac{1}{r^2} \right) \left( \frac{1}{r^{p + \delta_\ell}} \right)^\ell |\Omega|\\
\overset{\eqref{eq:DefinitionOfCRP} \eqref{eq:Easy}}&{\leq} \left( 1 - \frac{1}{r^2} \right) \left( 2 C_{r,p} \left( 1 - \frac{1}{r^2} \right) (t_{\ell + 1} - t_{\ell}) + \left( \frac{1}{r^{p + \delta_{\ell}}} \right)^{\ell} \right) |\Omega| \\
&\leq \left( 1 - \frac{1}{r^2} \right) \left( 2 C_{r,p} (t_{\ell + 1} - t_{\ell}) + \left( \frac{1}{r^{p + \delta_{\ell}}} \right)^{\ell} \right) |\Omega| \overset{\eqref{eq:Difficult}}{\leq} \left( 1 - \frac{1}{r^2} \right) \left(\frac{1}{r^{p + \delta_{\ell + 1}}} \right)^\ell |\Omega|; \\
\mathcal{Q}^{\text{Step $(\ell + 1)$}}_{I+\ell} \overset{\ref{item:CIS_J3}}&{\geq} \left( 1 - \frac{1}{r^p} \right) |\{ Dv_{I+\ell} \in \mathcal{V}_{I+\ell} \}| \overset{\ref{item:CIS_G4}}{\geq} \left( 1 - \frac{1}{r^p} \right) \left( \frac{1}{r^2} \right)^\ell |\Omega|; \\
\end{align*}
Finally, we also compute upper and lower bounds on $\widetilde{\mathcal{Q}}^{\text{Step $(\ell+1)$}}$:
\begin{align*}
\widetilde{\mathcal{Q}}^{\text{Step $(\ell+1)$}} \overset{\ref{item:CIS_H5}, \ref{item:CIS_I5}, \ref{item:CIS_J4}}&{\leq} \sum_{k = 0}^{\ell-1} \left( 1 - \frac{t_{\ell}}{t_{\ell + 1}} \right) \left( \frac{1}{r^{p + 2 \delta}} \right)^{\ell + 1 - k}  |\{ Dv_{I+\ell} \in \mathcal{W}_{I+k, t_{\ell}}^1 \}| \\
&\qquad + \sum_{k = 0}^{\ell-1}  \left( 1 - \frac{t_{\ell}}{t_{\ell + 1}} \right) \left(\frac{1}{r^{p + 2 \delta}} \right)^{\ell - k} |\{ Dv_{I+\ell} \in \mathcal{W}_{I+k, t_{\ell}}^2 \}| + \frac{1}{r^{p + 2 \delta}} |\{ Dv_{I+\ell} \in \mathcal{V}_{I+\ell} \}| \\
&\leq \left( 1 - \frac{t_{\ell}}{t_{\ell + 1}} \right) \sum_{k = 0}^{\ell-1} \left( \frac{1}{r^{p + \delta_1}} \right)^{\ell - k} |\{ Dv_{I+\ell} \in \mathcal{W}_{I+k, t_{\ell}}^1 \cup \mathcal{W}_{I+k, t_{\ell}}^2 \}| + \frac{1}{r^{p + \delta_1}} |\{ Dv_{I+\ell} \in \mathcal{V}_{I+\ell} \}| \\
\overset{\ref{item:CIS_G3}, \ref{item:CIS_G4}}&{\leq} \left( 1 - \frac{1}{r^2} \right) \left(1 - \frac{t_{\ell}}{t_{\ell+1}} \right) \left(\frac{1}{r^{p + \delta_1}}\right)^{\ell} \sum_{k = 0}^{\ell-1} \left(\frac{r^{p + \delta_1}}{r^{p + \delta_\ell}} \right)^{k} |\Omega| + \left(\frac{1}{r^{p + \delta_\ell}}\right)^{\ell+1} |\Omega| \\
&\leq \left( 1 - \frac{1}{r^2} \right) \left(1 - \frac{t_{\ell}}{t_{\ell+1}} \right) \ell \left(\frac{r^{\delta_1 - \delta_{\ell}}}{r^{p + \delta_1}}\right)^{\ell-1} |\Omega| + \left(\frac{1}{r^{p + \delta_\ell}}\right)^{\ell+1} |\Omega| \\
\overset{\eqref{eq:DefinitionOfCRP}, \eqref{eq:Easy}}&{\leq} \left( 2 C_{r,p} (t_{\ell + 1} - t_{\ell}) + \left(\frac{1}{r^{p + \delta_\ell}}\right)^{\ell+1} \right) |\Omega|\overset{\eqref{eq:Difficult}}{\leq} \left( \frac{1}{r^{p + \delta_{\ell+1}}} \right)^{\ell+1} |\Omega|; \\
\widetilde{\mathcal{Q}}^{\text{Step $(\ell+1)$}} \overset{\ref{item:CIS_J4}}&{\geq} \frac{1}{r^2} |\{ Dv_{I+\ell} \in \mathcal{V}_{I+\ell} \}| \overset{\ref{item:CIS_G4}}{\geq} \left( \frac{1}{r^2} \right)^{\ell+1} |\Omega|.
\end{align*}
Select $0 < \eps_{I + \ell + 1} < \min(\eps_{I + \ell}, 2^{-(I + \ell + 1)})$ so that ${\| D v_{I + \ell + 1} \ast \rho_{\eps_{I + \ell + 1}} - D v_{I + \ell + 1} \|_{L^1(\Omega; \R^{2 \times 2})} < 2^{-(I + \ell + 1)}}$. 
Let $\mathcal{F}_{I + \ell + 1}$ be a collection of open convex mutually disjoint sets covering $\Omega$ up to a set of measure 0 such that $v_{I + \ell + 1} |_{U}$ is affine for each $U \in \mathcal{F}_{I + \ell + 1}$ and $\sup_{U \in \mathcal{F}_{I + \ell + 1}} \diam{U} < 2^{-(I + \ell + 1)}$.
This ends the $(\ell + 1)$-th step. 
Therefore, by induction, we have built a sequence $\{ v_{I + \ell} \}_{\ell = 1}^{\infty}$ and $\{ \eps_{I + \ell} \}_{\ell = 1}^{\infty}$ such that for all $\ell \geq 1$:
\begin{align}
&\| v_{I+\ell} - v_{I + \ell -1 }\|_{C^{\alpha}(\overline{\Omega}; \R^2)} < \frac{\eps_{I + \ell - 1}}{2^{I + \ell - 1}}; \label{eq:HolderIncrement}\\
&Dv_{I + \ell}(x) \in \bigcup_{j = 0}^{\ell - 1} \mathcal{W}_{I+j,t_\ell}^1 \cup \bigcup_{j = 0}^{\ell - 1} \mathcal{W}_{I+j,t_\ell}^2 \cup \mathcal{V}_{I+\ell} \text{ for a.e. } x \in \Omega; \label{eq:DiffInclusionOfSeq} \\
&\left( 1 - \frac{1}{r^p} \right) \left( \frac{1}{r^2} \right)^j \frac{t_{j+1}}{t_{\ell}} |\Omega| \leq |\{ Dv_{I + \ell} \in \mathcal{W}_{I+j, t_\ell}^1 \cup \mathcal{W}_{I+j, t_\ell}^2 \}| \leq \left( 1 - \frac{1}{r^2} \right) \left( \frac{1}{r^{p + \delta_\ell}} \right)^j |\Omega|, \, \forall 0 \leq j < \ell; \label{eq:UpperAndLowerBoundsOfMassInW}\\
&\left( \frac{1}{r^2} \right)^\ell |\Omega| \leq |\{ Dv_{I + \ell} \in \mathcal{V}_{I+\ell} \}| \leq \left( \frac{1}{r^{p + \delta_\ell}} \right)^\ell |\Omega| \label{eq:UpperAndLowerBoundOnError};\\
&\| D v_{I + \ell} \ast \rho_{\eps_{I + \ell}} - D v_{I + \ell} \|_{L^1(\Omega; \R^{2 \times 2})} < 2^{-(I + \ell)}; \label{eq:ConvGradIneq}
\end{align}
and there is a collection $\mathcal{F}_{I+\ell}$ of open convex mutually disjoint sets covering $\Omega$ up to set of measure 0 such that $v_{I+\ell} |_U$ is affine for all $U \in \mathcal{F}_{I+\ell}$ and $\sup_{U \in \mathcal{F}_{I+\ell}} \diam{U} < 2^{-(I+\ell)}$.

Due to computations that we will carry out in Subsection~\ref{subsec:UniformSobolevBound}, we find that there is some $p > 1$ such that $\sup_{\ell \geq 1} \| v_{I + \ell} \|_{W^{1,p}(\Omega; \R^2)} < \infty$. Therefore there exists a function $v \in W^{1,p}(\Omega; \R^2)$ such that $v_{I + \ell}$ converges weakly to $v$ in $W^{1,p}(\Omega; \R^2)$ (no need to pass to subsequences since $v \in C^{\alpha}(\overline{\Omega}; \R^2)$ and $v_{I + \ell} \to v$ in $C^{\alpha}(\overline{\Omega}; \R^2)$). In Subsection~\ref{subsec:L1Convergence}, we prove that $Dv_{I + \ell} \to Dv$ in $L^1(\Omega; \R^{2 \times 2})$ as $\ell \to \infty$. Then we show that $Dv(x) \in K_{f}$ for a.e. $x \in \Omega$ in Subsection~\ref{subsec:DifferentialIncl}. This proves, by virtue of the arguments from Subsection~\ref{subsec:EulerLagrangeAsDifferentialInclusion}, that the first component of $v$, denoted by $v^1$, solves equation \eqref{eq:Intro:EL}. In Subsection~\ref{subsec:L2Inf}, we show that for any ball $B \subset \Omega$ we have
\[
 \int_{B} |\partial_1 v^1|^2 \, dx = + \infty.
\]
In Subsection~\ref{subsec:InfinitelyManySolutions}, we explain why our scheme yields infinitely many such solutions to \eqref{eq:Intro:EL} with the same boundary datum.
This ends the proof of Theorem~\ref{thm:MainFunctionals}.

\subsection{Proof of existence of the sequences $\{ \delta_{\ell} \}_{\ell = 1}^{\infty}$ and $\{ t_{\ell} \}_{\ell = 1}^{\infty}$} \label{subsec:tlseq}
Define the sequence $\{ \delta_{\ell} \}_{\ell = 1}^{\infty}$ by
\[
 \delta_{\ell} \coloneqq \delta + \frac{1}{\ell^2} \delta = \left( 1 + \frac{1}{\ell^2} \right) \delta.
\]
By standard computations, $\delta_{\ell} - \delta_{\ell + 1} \geq \frac{2}{(\ell+1)^3} \delta$. This implies that
\begin{equation}\label{eq:SeqIncrementBoundFromBelow}
 (\delta_{\ell} - \delta_{\ell + 1}) (\ell + 1)^3 \geq 2 \delta.
\end{equation}
In addition,
\[
 (\ell + 1)^2 \max \left\{ \left( \dfrac{1}{r^{p + \delta_{\ell+1}}} \right)^{(\ell + 1)}, \left( \dfrac{1}{r^{p + \delta_{\ell}}} \right)^{(\ell + 1)} \right\} \to 0 \quad \text{ as $\ell \to \infty$.}
\]
This implies that there exists a constant $\widetilde{C} \geq 1$ such that
\begin{equation}\label{eq:QuantityBoundedByTheConstantM}
 (\ell + 1)^2 \max \left\{ \left( \dfrac{1}{r^{p + \delta_{\ell+1}}} \right)^{(\ell + 1)}, \left( \dfrac{1}{r^{p + \delta_{\ell}}} \right)^{(\ell + 1)} \right\} \leq \widetilde{C} \quad \forall \ell \geq 1.
\end{equation}
Define the sequence $\{ t_{\ell} \}_{\ell = 1}^{\infty}$ by
\[
 t_{\ell} \coloneqq 1 - \frac{1}{r^{N + I + 2p(\ell + 1)^3} (\ell + 1)}.
\]
where $N \geq \log_r (2 \widetilde{C} C_{r,p})$ is selected sufficiently large so that \eqref{eq:Easy} holds. The fact that \eqref{eq:SeqLowerBound} holds is immediate. It remains to prove \eqref{eq:Difficult}. We have
\begin{equation}\label{eq:BoundIncrementOfTSeq}
 t_{\ell + 1} - t_{\ell} \leq \frac{1}{r^{N + I + 2p(\ell + 1)^3} (\ell + 1)}
\end{equation}
and
\begin{align*}
 \left( \frac{1}{r^{p + \delta_{\ell + 1}}} \right)^{(\ell + 1)^3} - \left( \frac{1}{r^{p + \delta_{\ell}}} \right)^{(\ell + 1)^3} &= \frac{(r^{p + \delta_{\ell}})^{(\ell + 1)^3} - (r^{p + \delta_{\ell + 1}})^{(\ell + 1)^3} }{(r^{2p + \delta_{\ell} + \delta_{\ell + 1}})^{(\ell + 1)^{3}}} = \frac{1}{r^{p(\ell + 1)^3}} \frac{r^{(\delta_{\ell} - \delta_{\ell + 1})(\ell + 1)^3} - 1}{r^{\delta_{\ell} (\ell + 1)^3 }} \\
 &\geq \frac{1}{r^{2 p (\ell + 1)^3}} (r^{(\delta_{\ell} - \delta_{\ell + 1})(\ell + 1)^3} - 1) \overset{\eqref{eq:SeqIncrementBoundFromBelow}}{\geq} \frac{1}{r^{2 p (\ell + 1)^3}} (r^{2 \delta} - 1) \\
 \overset{\eqref{eq:AboutIToThePowerTwo}}&{\geq} \frac{1}{r^{I + 2 p (\ell + 1)^3}} = r^N (\ell + 1) \frac{1}{r^{N + I + 2 p (\ell + 1)^3}(\ell + 1)} \\
 \overset{\eqref{eq:BoundIncrementOfTSeq}}&{\geq} 2 \widetilde{C} C_{r,p} (t_{\ell + 1} - t_{\ell}) (\ell + 1). \\
\end{align*}
Therefore, 
\begin{equation}\label{eq:SomeBoundForLPulsOneThatWillBeUsed}
\begin{split}
&2 C_{r,p} (t_{\ell + 1} - t_{\ell}) (\ell + 1) \\
&\quad \leq \frac{1}{\widetilde{C}} \left( \left( \frac{1}{r^{p + \delta_{\ell + 1}}} \right)^{(\ell + 1)^3} - \left( \frac{1}{r^{p + \delta_{\ell}}} \right)^{(\ell + 1)^3} \right) \\
&\quad \leq \frac{1}{\widetilde{C}} (\ell + 1)^2 \max \left\{ \left( \frac{1}{r^{p + \delta_{\ell + 1}}} \right)^{(\ell + 1)} , \left( \frac{1}{r^{p + \delta_{\ell}}} \right)^{(\ell + 1)} \right\}^{(\ell + 1)^2 - 1} \left( \left( \frac{1}{r^{p + \delta_{\ell + 1}}} \right)^{(\ell + 1)} - \left( \frac{1}{r^{p + \delta_{\ell}}} \right)^{(\ell + 1)} \right) \\
&\quad \leq \frac{1}{\widetilde{C}} (\ell + 1)^2 \max \left\{ \left( \frac{1}{r^{p + \delta_{\ell + 1}}} \right)^{(\ell + 1)} , \left( \frac{1}{r^{p + \delta_{\ell}}} \right)^{(\ell + 1)} \right\} \left( \left( \frac{1}{r^{p + \delta_{\ell + 1}}} \right)^{(\ell + 1)} - \left( \frac{1}{r^{p + \delta_{\ell}}} \right)^{(\ell + 1)} \right) \\
&\quad \overset{\eqref{eq:QuantityBoundedByTheConstantM}}{\leq} \left( \frac{1}{r^{p + \delta_{\ell + 1}}} \right)^{(\ell + 1)} - \left( \frac{1}{r^{p + \delta_{\ell}}} \right)^{(\ell + 1)}
\end{split}
\end{equation}
for all $\ell \geq 1$, where we used the fact that
\begin{equation}\label{eq:ForRealNumber-p-ContractionObservation}
 |x^q - y^q| \leq q \max \{ |x| , |y| \}^{q-1} |x-y| \quad \forall x,y \in [0,1], \, \forall q \geq 1
\end{equation}
and
\begin{equation}\label{eq:AMaximumThatIsBoundedByOne}
 \max \left\{ \left( \frac{1}{r^{p + \delta_{\ell + 1}}} \right)^{(\ell + 1)} , \left( \frac{1}{r^{p + \delta_{\ell}}} \right)^{(\ell + 1)} \right\} \leq 1 \quad \forall \ell \geq 1.
\end{equation}
Using \eqref{eq:ForRealNumber-p-ContractionObservation} and \eqref{eq:AMaximumThatIsBoundedByOne} once again, we find that for all $j = 1, \ldots, \ell + 1$,
\begin{align*}
 2 C_{r,p} (t_{\ell + 1} - t_{\ell}) \overset{\eqref{eq:SomeBoundForLPulsOneThatWillBeUsed}}&{\leq} \frac{1}{\ell + 1} \left( \left( \frac{1}{r^{p + \delta_{\ell + 1}}} \right)^{(\ell + 1)} - \left( \frac{1}{r^{p + \delta_{\ell}}} \right)^{(\ell + 1)} \right) \\
 &= \frac{1}{\ell + 1} \left( \left( \frac{1}{r^{p + \delta_{\ell + 1}}} \right)^{j \frac{(\ell + 1)}{j}} - \left( \frac{1}{r^{p + \delta_{\ell}}} \right)^{j \frac{(\ell + 1)}{j}} \right) \\
 \overset{\eqref{eq:ForRealNumber-p-ContractionObservation}}&{\leq} \frac{1}{\ell + 1} \frac{\ell + 1}{j} \left( \left( \frac{1}{r^{p + \delta_{\ell + 1}}} \right)^{j} - \left( \frac{1}{r^{p + \delta_{\ell}}} \right)^{j} \right)\leq \left( \frac{1}{r^{p + \delta_{\ell + 1}}} \right)^{j} - \left( \frac{1}{r^{p + \delta_{\ell}}} \right)^{j}.\\
\end{align*}
This proves \eqref{eq:Difficult}.

\subsection{Proof that $\sup_{\ell \geq 1} \| v_{I + \ell} \|_{W^{1,p}(\Omega; \R^2)} < \infty$}\label{subsec:UniformSobolevBound}
First we notice that due to \eqref{eq:ExpontialBoundInTheSets}
\[
 \sup_{\ell \geq 1} \| \partial_2 v_{I + \ell}^1 \|_{L^{\infty}(\Omega; \R^{2 \times 2})} , \quad \sup_{\ell \geq 1} \| \partial_1 v_{I + \ell}^2 \|_{L^{\infty}(\Omega; \R^{2 \times 2})}< 1. 
\] 
In addition, we have (with $C$ given by \eqref{eq:ExpontialBoundInTheSets})
\begin{align*}
 \int_{\Omega} &|\partial_1 v_{I + \ell}^1|^p + |\partial_2 v_{I + \ell}^2|^p \, dx = \sum_{j = 0}^{\ell - 1} \int_{\{ Dv_{I + \ell} \in \mathcal{W}_{I+j, t_\ell}^1 \} \cup \{ Dv_{I + \ell} \in \mathcal{W}_{I+j, t_\ell}^2 \}} |\partial_1 v_{I + \ell}^1|^p + |\partial_2 v_{I + \ell}^2|^p \, dx \\
 &\qquad + \int_{\{ Dv_{I + \ell} \in \mathcal{V}_{I + \ell} \}} |\partial_1 v_{I + \ell}^1|^p + |\partial_2 v_{I + \ell}^2|^p \, dx \\
 \overset{\eqref{eq:ExpontialBoundInTheSets}}&{\leq} C r^{p I} \sum_{j = 0}^{\ell - 1} r^{pj} |\{ Dv_{I + \ell} \in \mathcal{W}_{I+j, t_\ell}^1 \cup \mathcal{W}_{I+j, t_\ell}^2 \}| + Cr^{p (I + \ell)} |\{ Dv_{I + \ell} \in \mathcal{V}_{I+\ell} \}| \\
 \overset{\eqref{eq:UpperAndLowerBoundsOfMassInW}, \eqref{eq:UpperAndLowerBoundOnError}}&{\leq} C r^{p I} |\Omega| \sum_{j = 0}^{\ell - 1} r^{pj} \left( 1 - \frac{1}{r^2} \right) \left( \frac{1}{r^{p + \delta_\ell}} \right)^j + Cr^{p ( I + \ell)} |\Omega| \left( \frac{1}{r^{p + \delta_\ell}} \right)^\ell \leq C r^{p I} |\Omega| \sum_{j = 0}^{\infty} \left( \frac{r^p}{r^{p + \delta}} \right)^j < \infty.
\end{align*}
Therefore, since the last term does not depend on $\ell$, we conclude that $\sup_{\ell \geq 1} \| v_{I + \ell} \|_{W^{1,p}(\Omega; \R^2)} < \infty$.

\subsection{Proof that $Dv_{I + \ell} \to Dv$ in $L^1(\Omega; \R^{2 \times 2})$ as $\ell \to \infty$}\label{subsec:L1Convergence}
We have, by \eqref{eq:ConvGradIneq},
\begin{equation}\label{eq:L1ConvIneq}
\begin{split}
 \| Dv_{I + \ell} - Dv \|_{L^1(\Omega; \R^{2 \times 2})} &\leq \underbrace{\| Dv_{I + \ell} - Dv_{I + \ell} \ast \rho_{\eps_{I + \ell}} \|_{L^1(\Omega; \R^{2 \times 2})}}_{\leq 2^{-(I + \ell)}} + \| Dv_{I + \ell} \ast \rho_{\eps_{I + \ell}} - Dv \ast \rho_{\eps_{I + \ell}} \|_{L^1(\Omega; \R^{2 \times 2})} \\
 &\qquad + \underbrace{\| Dv \ast \rho_{\eps_{I + \ell}} - Dv \|_{L^1(\Omega; \R^{2 \times 2})}.}_{\to 0 \text{ as } \ell \to \infty}
\end{split}
\end{equation}
The first and third term converge to 0 as $\ell \to \infty$. Therefore, it suffices to show that the second term converges to 0 as $\ell \to \infty$. Inside $\Omega$ (recall the conventions fixed in Section~\ref{sec:Notation} and that $v_{I + \ell}$ has affine boundary datum), we have $(Dv_{I + \ell} - Dv)\ast \rho_{\eps_{I + \ell}} = (v_{I + \ell} - v) \ast D\rho_{\eps_{I + \ell}}$. We have
\begin{align*}
 \| (v_{I + \ell} - v) \ast D\rho_{\eps_{I + \ell}} \|_{L^1(\Omega; \R^{2 \times 2})} \leq \frac{\| D \rho \|_{L^1(B_1; \R^{2})}}{\eps_{I + \ell}} \| v_{I + \ell} - v \|_{L^1(\Omega; \R^{2})} \leq \frac{\| D \rho \|_{L^1(B_1; \R^{2})} |\Omega|}{\eps_{I + \ell}} \| v_{I + \ell} - v \|_{L^{\infty}(\Omega; \R^{2})}.
\end{align*}
In addition,
\begin{equation*}
 \| v_{I + \ell} - v \|_{L^{\infty}(\Omega; \R^{2})} \leq \sum_{j = I + \ell}^{\infty} \| v_{j+1} - v_{j} \|_{L^{\infty}(\Omega; \R^{2})} \overset{\eqref{eq:HolderIncrement}}{\leq} \sum_{j = I + \ell}^{\infty} \frac{\eps_{I+j}}{2^{I+j}} \leq \frac{\eps_{I+\ell}}{2^{I+\ell-1}}.
\end{equation*}
Thus the second term in \eqref{eq:L1ConvIneq} converges to 0 as $\ell \to \infty$. From this, we deduce that $Dv_{I + \ell} \to Dv$ in $L^1(\Omega; \R^{2 \times 2})$ as $\ell \to \infty$.

\subsection{Proof that $Dv(x) \in K_{f}$ for a.e. $x \in \Omega$}\label{subsec:DifferentialIncl}
We begin by proving that
\begin{equation*}
 Dv(x) \in \underbrace{\bigcup_{j = 0}^{\infty} \overline{\mathcal{U}}_{I+j}^1 \cup \bigcup_{j = 0}^{\infty} \overline{\mathcal{U}}_{I+j}^2 \cup \bigcup_{j = 0}^{\infty} \overline{\mathcal{V}}_{I+j}}_{\qquad =: \tilde{K}_{f}} \quad \text{for a.e. } x \in \Omega.
\end{equation*}
From \eqref{eq:SeqLowerBound}, we find that
\begin{equation}\label{eq:InApproxDist}
\dist{X}{\tilde{K}_{f}} \leq \frac{1}{2^{\ell}} \quad \forall X \in \bigcup_{j = 0}^{\ell - 1} \mathcal{W}_{I+j, t_{\ell}}^1 \cup \bigcup_{j = 0}^{\ell - 1} \mathcal{W}_{I+j, t_{\ell}}^2 \cup \mathcal{V}_{I+\ell}.
\end{equation}
The fact that $Dv_{I + \ell}(x) \in \bigcup_{j = 0}^{\ell - 1} \mathcal{W}_{I+j, t_{\ell}}^1 \cup \bigcup_{j = 0}^{\ell - 1} \mathcal{W}_{I+j, t_{\ell}}^2 \cup \mathcal{V}_{I+\ell}$ for a.e. $x \in \Omega$, combined with \eqref{eq:InApproxDist} and the fact that $Dv_{I+\ell} \to Dv$ in $L^1(\Omega; \R^{2 \times 2})$ allows us to conclude that $Dv(x) \in \tilde{K}_{f}$ for a.e. $x \in \Omega$. It remains to prove that
\begin{equation}\label{eq:ErrorsCancel}
 \left| \left\{ Dv(x) \in \bigcup_{j = 0}^{\infty} \overline{\mathcal{V}}_{I + j} \right\} \right| = 0.
\end{equation}
First we point out that since $Dv_{I+\ell} \to Dv$ in $L^1(\Omega; \R^{2 \times 2})$, $Dv_{I+\ell} \to Dv$ in measure. Notice that
\begin{align*}
 \left\{ Dv(x) \in \bigcup_{j = 0}^{\infty} \overline{\mathcal{V}}_{I + j} \right\} &\subset \left\{ Dv(x) \in \bigcup_{j = 0}^{\infty} \overline{\mathcal{V}}_{I + j}, |Dv_{I+\ell} - Dv| > 1 \right\} \\
 & \qquad \cup \left\{ Dv(x) \in \bigcup_{j = 0}^{\infty} \overline{\mathcal{V}}_{I + j}, |Dv_{I+\ell} - Dv| \leq 1 \right\} \\
 \overset{\eqref{eq:AboutDistancesBetweenSets}}&{\subset}  \left\{ |Dv_{I+\ell} - Dv| > 1 \right\} \cup \left\{ Dv_{I + \ell} \in \mathcal{V}_{I + \ell}, \right\}.
\end{align*}
Therefore
\begin{align*}
 \left| \left\{ Dv(x) \in \bigcup_{j = 0}^{\infty} \overline{\mathcal{V}}_{I + j} \right\} \right| &\leq \left| \left\{ |Dv_{I+\ell} - Dv| > 1 \right\} \right| + \left| \left\{ Dv_{I + \ell} \in \mathcal{V}_{I + \ell}, \right\} \right| \\
 \overset{\eqref{eq:UpperAndLowerBoundOnError}}&{\leq} \left| \left\{ |Dv_{I+\ell} - Dv| > 1 \right\} \right| + \left( \frac{1}{r^{p + \delta}} \right)^{\ell} |\Omega|.
\end{align*}
Since $Dv_{I+\ell} \to Dv$ in measure, the first term converges to 0 as $\ell \to \infty$.
Since the second term converges to 0 as well, we find \eqref{eq:ErrorsCancel}, as wished. 
From \eqref{eq:ErrorsCancel} we conclude that 
\begin{equation*}
 Dv(x) \in \bigcup_{j = 0}^{\infty} \overline{\mathcal{U}}_{I+j}^1 \cup \bigcup_{j = 0}^{\infty} \overline{\mathcal{U}}_{I+j}^2 \subset K_{f} \quad \text{for a.e. } x \in \Omega.
\end{equation*}

\subsection{Proof that $\int_{B} |\partial_1 v^1|^2 \, dx = + \infty$ for any ball $B \subset \Omega$}\label{subsec:L2Inf}
Let $B$ be an open ball of arbitrary radius. Let $B^{\prime} \subset B$ be an arbitrary ball which is a subset of $B$. Since $\sup_{U \in \mathcal{F}_{I+\ell}} \operatorname{diam}(U) < 2^{- (I + \ell)}$ for all $\ell \geq 1$, for $\ell$ large enough there is an open set $U \in \mathcal{F}_{I+\ell}$ so that $U \subset B^{\prime}$. By construction, there is a subset $\tilde{U} \in \mathcal{F}_{I + \ell + 1}$, $\tilde{U} \subset U$ so that $v_{I + \ell + 1}|_{\tilde{U}}$ is affine and $Dv_{I + \ell + 1}|_{\tilde{U}} \in \mathcal{V}_{I + \ell + 1}$. 
By construction, there exist $\eps_{\tilde{U}} > 0$ and $\eps^{\prime}_{\tilde{U}} > 0$ such that
\begin{align*}
|\{ Dv_{I + \ell + 3} |_{\tilde{U}} \in \mathcal{W}_{I + \ell + 1, t_{\ell + 3}}^1 \cup \mathcal{W}_{I + \ell + 1, t_{\ell + 3}}^2 \}| > \eps_{\tilde{U}} |\tilde{U}|; \\
|\{ Dv_{I + \ell + 3} |_{\tilde{U}} \in \mathcal{W}_{I + \ell + 2, t_{\ell + 3}}^1 \cup \mathcal{W}_{I + \ell + 2, t_{\ell + 3}}^2 \}| > \eps_{\tilde{U}}^{\prime} |\tilde{U}|.
\end{align*}
By \ref{item:CIS_H3} and \ref{item:CIS_I3}, we obtain 
\begin{align*}
|\{ Dv_{I + \ell + 4} |_{\tilde{U}} \in \mathcal{W}_{I + \ell + 1, t_{\ell + 4}}^1 \cup \mathcal{W}_{I + \ell + 1, t_{\ell + 4}}^2 \}| > \frac{t_{\ell + 3}}{t_{\ell + 4}} \eps_{\tilde{U}} |\tilde{U}|; \\
|\{ Dv_{I + \ell + 4} |_{\tilde{U}} \in \mathcal{W}_{I + \ell + 2, t_{\ell + 4}}^1 \cup \mathcal{W}_{I + \ell + 2, t_{\ell + 4}}^2 \}| > \frac{t_{\ell + 3}}{t_{\ell + 4}}\eps_{\tilde{U}}^{\prime} |\tilde{U}|.
\end{align*}
Letting $s \geq 5$ be an arbitrary integer and applying \ref{item:CIS_H3} and \ref{item:CIS_I3} $(s - 4)$ more times yield, thanks to \eqref{eq:Easy}
\begin{align*}
|\{ Dv_{I + \ell + s} |_{\tilde{U}} \in \mathcal{W}_{I + \ell + 1, t_{\ell + s}}^1 \cup \mathcal{W}_{I + \ell + 1, t_{\ell + s}}^2 \}| > \frac{t_{\ell + 3}}{t_{\ell + s}} \eps_{\tilde{U}} |\tilde{U}| > \frac{9}{10} \eps_{\tilde{U}} |\tilde{U}|; \\
|\{ Dv_{I + \ell + s} |_{\tilde{U}} \in \mathcal{W}_{I + \ell + 2, t_{\ell + s}}^1 \cup \mathcal{W}_{I + \ell + 2, t_{\ell + s}}^2 \}| > \frac{t_{\ell + 3}}{t_{\ell + s}}\eps_{\tilde{U}}^{\prime} |\tilde{U}| > \frac{9}{10} \eps_{\tilde{U}}^{\prime} |\tilde{U}|.
\end{align*}
Due to the fact that $Dv_{I + \ell + s} \to Dv$ in $L^1(\Omega; \R^{2 \times 2})$ as $s \to \infty$, this implies that
\begin{align*}
|\{ Dv |_{\tilde{U}} \in \overline{\mathcal{U}}_{I + \ell + 1}^1 \cup \overline{\mathcal{U}}_{I + \ell + 1}^2 \}| \geq \frac{9}{10} \eps_{\tilde{U}} |\tilde{U}| \quad \text{and} \quad 
|\{ Dv |_{\tilde{U}} \in \overline{\mathcal{U}}_{I + \ell + 2}^1 \cup \overline{\mathcal{U}}_{I + \ell + 2}^2  \}| \geq \frac{9}{10} \eps_{\tilde{U}}^{\prime} |\tilde{U}|.
\end{align*}
Thus, by assuming $\ell$ sufficiently large so that
\[
 \dist{ \overline{\mathcal{U}}_{I + \ell + 1}^1 \cup \overline{\mathcal{U}}_{I + \ell + 1}^2 }{ \overline{\mathcal{U}}_{I + \ell + 2}^1 \cup \overline{\mathcal{U}}_{I + \ell + 2}^2 } > 1
\]
we obtain
\begin{equation*}
 \operatorname{ess\,sup}_{x,y \in B^{\prime}} |\partial_1 v^1(x) - \partial_1 v^1(y)| > 1.
\end{equation*}
Since this last inequality holds for all $B^{\prime} \subset B$, $v^1$ is not $C^1$ on $B$. Since $f$ is smooth, due to the results by De Giorgi and Nash \cite{EDG, JN58} $v^1 |_B \not\in W^{1,2}(B)$. Indeed, otherwise, if $v^1 |_B \in W^{1,2}(B)$, then by these results we deduce that $v^1 \in C^{\infty}(B)$. This is a contradiction.

\subsection{Why does this convex integration scheme yield infinitely many solutions with identical boundary datum?}\label{subsec:InfinitelyManySolutions}
To see that the convex integration scheme yields infinitely many solutions, all having the same boundary datum, we first observe that with $v$ given by the scheme above
\[
 \| v - v_I \|_{C^{\alpha}(\overline{\Omega}; \R^2)} \leq \sum_{j > I} \| v_j - v_{j-1} \|_{C^{\alpha}(\overline{\Omega}; \R^2)} \overset{\eqref{eq:HolderIncrement}}{\leq} \sum_{j > I} \frac{\eps_j}{2^j} < \eps_I \sum_{j > I} \frac{1}{2^j} < \eps_I < \eps.
\]
Keep in mind that $\eps > 0$ can be chosen arbitrarily at the beginning of the convex integration scheme. Assume for a contradiction that it is only possible to find finitely many $v^{(i)}$ ($i = 1 , \ldots, M$) such that $Dv^{(i)} \in K_{f}$ a.e. Then select (recall that $v_I(x) = Mx$)
\[
 \tilde{\eps} \coloneqq \frac{\min_{i = 1, \ldots, M} \| v^{(i)} - v_I \|_{C^{\alpha}(\overline{\Omega}; \R^2)}}{2}.
\]
By the convex integration scheme, we can construct a function $\tilde{v} \in W^{1,p}(\Omega; \R^2) \cap C^{\alpha}(\overline{\Omega}; \R^2)$ such that $\| \tilde{v} - v_I \|_{C^{\alpha}(\overline{\Omega}; \R^2)} < \tilde{\eps}$ and $D \tilde{v} \in K_{f}$ a.e.. This is a contradiction. This proves that there are infinitely many solutions to \eqref{eq:Intro:EL}, all with identical boundary datum.

\bibliographystyle{plain}
\bibliography{ReferencesNonenergetic}

\end{document}